\input amstex\documentstyle {amsppt}  
\pagewidth{12.5 cm}\pageheight{19 cm}\magnification\magstep1
\topmatter
\title Character sheaves on disconnected groups, VI\endtitle
\author G. Lusztig\endauthor
\address Department of Mathematics, M.I.T., Cambridge, MA 02139\endaddress
\thanks Supported in part by the National Science Foundation. The author thanks the 
Institut des Hautes \'Etudes Scientifiques for its hospitality during a visit (may 
2004) when part of the work on this paper was done.\endthanks
\endtopmatter   
\document
\define\ufs{\un{\fs}}

\define\tius{\ti{\us}}

\define\dw{\dot w}
\define\dt{\dot t}
\define\ds{\dot s}
\define\dy{\dot y}
\define\bcl{\bar{\cl}}
\define\us{\un s}
\define\uD{\un{D}}

\define\dts{\dot{\ts}}

\define\bn{\bar n}

\define\bF{\bar F}

\define\da{\dagger}

\define\dsv{\dashv}

\define\po{\text{\rm pos}}

\define\si{\sim}

\define\sqc{\sqcup}

\define\qua{\quad}
\define\hA{\hat A}

\define\tcl{\ti\cl}

\define\til{\ti\l}
\define\tiz{\ti\z}
\define\tid{\ti\d}

\define\bK{\bar K}

\define\bZ{\bar Z}

\define\bpi{\bar\p}

\define\lb{\linebreak}

\define\op{\oplus}

\define\em{\emptyset}
\define\imp{\implies}
\define\ra{\rangle}
\define\n{\notin}
\define\iy{\infty}
\define\m{\mapsto}
\define\do{\dots}
\define\la{\langle}
\define\bsl{\backslash}

\define\lra{\leftrightarrow}

\define\sub{\subset}

\define\T{\times}
\define\ti{\tilde}
\define\nl{\newline}
\redefine\i{^{-1}}
\define\fra{\frac}
\define\un{\underline}
\define\ov{\overline}
\define\ot{\otimes}
\define\bbq{\bar{\QQ}_l}

\define\Ad{\text{\rm Ad}}
\define\Hom{\text{\rm Hom}}

\define\Aut{\text{\rm Aut}}

\define\ind{\text{\rm ind}}

\define\res{\text{\rm res}}

\define\supp{\text{\rm supp}}

\define\bst{\bigstar}

\define\a{\alpha}
\redefine\b{\beta}
\redefine\c{\chi}
\define\g{\gamma}
\redefine\d{\delta}
\define\e{\epsilon}

\define\io{\iota}
\redefine\o{\omega}
\define\p{\pi}
\define\ph{\phi}
\define\ps{\psi}
\define\r{\rho}

\redefine\t{\tau}

\define\k{\kappa}
\redefine\l{\lambda}
\define\z{\zeta}
\define\x{\xi}

\redefine\G{\Gamma}
\redefine\D{\Delta}

\define\dd{\bold d}

\define\kk{\bold k}

\redefine\ss{\bold s}
\redefine\tt{\bold t}
\define\uu{\bold u}

\define\ww{\bold w}

\define\yy{\bold y}

\redefine\AA{\bold A}

\define\FF{\bold F}

\define\II{\bold I}

\define\NN{\bold N}

\define\QQ{\bold Q}

\define\TT{\bold T}

\define\WW{\bold W}
\define\ZZ{\bold Z}

\define\ca{\Cal A}
\define\cb{\Cal B}

\define\cd{\Cal D}
\define\ce{\Cal E}
\define\cf{\Cal F}

\define\ch{\Cal H}

\define\cj{\Cal J}
\define\ck{\Cal K}
\define\cl{\Cal L}
\define\cm{\Cal M}

\define\cp{\Cal P}

\define\cs{\Cal S}
\define\ct{\Cal T}

\define\cz{\Cal Z}
\define\cx{\Cal X}

\define\fa{\frak a}

\define\fii{\frak i}

\define\fs{\frak s}

\define\fZ{\frak Z}

\define\te{\ti e}
\define\tf{\ti f}

\define\ts{\ti s}

\define\tG{\ti G}

\define\tZ{\ti Z}

\define\sh{\sharp}

\define\sps{\supset}

\define\bS{\bar S}

\define\bul{\bullet}

\define\che{\check}
\define\BBD{BBD}
\define\DE{D}
\define\CS{L3}
\define\AD{L9}
\define\PCS{L10}
\define\LM{L11}
\define\HA{L12}
\define\MS{MS}
\define\YO{Y}
\define\hfZ{\hat{\fZ}}

\define\cha{\che\a}

\head Introduction\endhead
Throughout this paper, $G$ denotes a fixed, not necessarily connected, reductive
algebraic group over an algebraically closed field $\kk$. This paper is a part of a
series \cite{\AD} which attempts to develop a theory of character sheaves on $G$. 

In Section 28 we define the character sheaves on a connected component of $G$ 
generalizing the definition in \cite{\CS, I,\S2}. In Section 29 we prove a
semisimplicity property of the restriction functor, generalizing one in 
\cite{\CS, I,\S3}. In Section 30 we show that any character sheaf is admissible, 
generalizing a result in \cite{\CS, I,\S4}. In Section 31 we show that the restriction
functor takes a character sheaf to a direct sum of character sheaves, generalizing a 
result in \cite{\CS, I,\S6}. 

We adhere to the notation of \cite{\AD} and \cite{\BBD}. Here is some additional 
notation. If $K\in\cd(X)$ and $A$ is a simple perverse sheaf on $X$ we write $A\dsv K$ 
instead of "$A$ is a subquotient of ${}^pH^i(K)$ for some $i\in\ZZ$." Let $\cm(X)$ be 
the subcategory of $\cd(X)$ whose objects are the perverse sheaves on $X$.

\head Contents\endhead
28. Definition of character sheaves. 

29. Restriction functor for character sheaves.

30. Admissibility of character sheaves.

31. Character sheaves and Hecke algebras.

\head 28. Definition of character sheaves\endhead
\subhead 28.1\endsubhead
Let $T$ be a torus. For any $n\in\NN^*_\kk$, let $\fs_n(T)$ be the category whose 
objects are the local systems of rank $1$ on $T$ that are equivariant for the 
transitive $T$-action $z:t\m z^nt$ on $T$; let $\fs(T)$ be the category whose objects 
are the local systems on $T$ that are in $\fs_n(T)$ for some $n$ as above.

If $f:T@>>>T'$ is a morphism of tori and $\cl'\in\fs(T')$ then $f^*\cl'\in\fs(T)$. The
set $\ufs(T)$ of isomorphism classes of objects in $\fs(T)$ is an abelian group for 
tensor product of local systems. Let $\cx=\Hom(T,\kk^*)$ (homomorphisms of algebraic 
groups). From the definitions we see that 

(a) {\it $\k\ot\ce\m\k^*\ce$ defines a group isomorphism 
$\cx\ot\ufs(\kk^*)@>\si>>\ufs(T)$.}
\nl
We show:

(b) {\it for $\cl\in\fs(T)$ there exists $\k\in\cx$ and $\ce\in\fs(\kk^*)$ such that
$\cl\cong\k^*\ce$.}
\nl
Indeed by (a) there exist $\k_i\in\cx,\ce_i\in\fs(\kk^*),(i\in[1,m])$ such that
$\cl\cong\ot_{i=1}^m\k_i^*\ce_i$. By 5.3 we have 
$\ufs(\kk^*)=\Hom(\mu_\iy(\kk^*),\bbq^*)\cong\QQ'/\ZZ$ where 
$\QQ'=\cup_{n\in\NN^*_\kk}\fra{1}{n}\ZZ\sub\QQ$. Hence we can find $\ce\in\fs(\kk^*)$,
$n_i\in\NN^*_\kk$ such that $\ce_i\cong\ce^{\ot n_i}$ for $i\in[1,m]$. Then 
$\cl\cong\ot_{i=1}^m\k_i^*\ce^{\ot n_i}\cong\k^*\ce$ where 
$\k=\prod_{i=1}^m\k_i^{n_i}$ and (b) follows.

For any $\t\in T$ define $h_\t:T@>>>T$ by $h_\t(t)=\t t$. We show:

(c) {\it if $\t\in T,\cl\in\fs(T)$, then $h_\t^*\cl\cong\cl$.}
\nl
Let $n$ be such that $\cl\in\fs_n(T)$. Then for any $z\in T$ we have
$h_{z^n}\cl\cong\cl$. We can find $z\in T$ such that $z^n=\t$. This proves (c).

\subhead 28.2\endsubhead
Let $\cl\in\fs(T)$, let $\k,\ce$ be as in 28.1(b) and let $n$ be the order of $\ce$ in
$\ufs(\kk^*)$. Then $\ce\in\fs_n(\kk^*)$. We show that the following two conditions 
for a morphism $f:T@>>>T$ of tori are equivalent:

(i) {\it $f^*\cl\cong\cl$;}

(ii) {\it there exists $\k_1\in\cx$ such that $\k\circ f=\k\k_1^n$.}
\nl
Condition (i) is equivalent to $f^*\k^*\ce\cong\k^*\ce$ that is,
$(\k\circ f)^*\ce\cong\k^*\ce$. Using the injectivity of the map 28.1(a) we see that 
this is equivalent to $(\k\circ f)\ot(n'/n)=\k\ot(n'/n)$ in $\cx\ot\QQ'/\ZZ$ (here 
$n'\in\ZZ,0<n'\le n$ and $n'/n$ is irreducible) which is clearly equivalent to 
condition (ii).

Assuming that (i),(ii) hold, we show:

(a) {\it $\cl$ is $T$-equivariant for the $T$-action $t_0:t\m f(t_0)tt_0\i$ on $T$.}
\nl
The map $\k:T@>>>\kk^*$ is compatible with the $T$-action (a) on $T$ and the $T$-action
$t_0:z\m\k_1(t_0)^nz$ on $\kk^*$. Hence to show that $\cl=\k^*\ce$ is $T$-equivariant 
it suffices to show that $\ce$ is $T$-equivariant. Since the $T$-action on $\kk^*$ 
comes via $\k_1$ from the $\kk^*$-action $z_0:z\m z_0^nz$ on $\kk^*$, it suffices to 
show that $\ce$ is $\kk^*$-equivariant. This holds since $\ce\in\fs_n(\kk^*)$.

\subhead 28.3\endsubhead
$G$ acts on $\prod_{B\in\cb}B/U_B$ by 
$$x:(g_BU_B)_{B\in\cb}\m(g'_BU_B)_{B\in\cb}$$
where $g'_{xBx\i}U_{xBx\i}=xg_Bx\i U_{xBx\i}$. Let 
$$\TT=(\prod_{B\in\cb}B/U_B)^{G^0}$$ 
(fixed point set of $G^0$). For any $B'\in\cb$ we define $f_{B'}:\TT@>\si>>B'/U_{B'}$ 
by $f_{B'}((g_BU_B)_{B\in\cb})=g_{B'}U_{B'}$. We use $f_{B'}$ to transport the 
algebraic group structure of $B'/U_{B'}$ to an algebraic group structure of $\TT$. This
structure is independent of the choice of $B'$. Thus $\TT$ is naturally a torus over 
$\kk$. The $G$ action on $\prod_{B\in\cb}B/U_B$ induces a $G/G^0$-action $D:t\m\uD(t)$ 
on $\TT$, respecting the algebraic group structure of $\TT$. We say that $\TT$ is the 
{\it canonical torus} of $G^0$.

For $w\in\WW$ (see 26.1) there is a unique isomorphism $\TT@>\si>>\TT$ (denoted again 
by $w$) such that for any $(B,B')\in\cb\T\cb$ with $\po(B,B')=w$ we have a commutative 
diagram
$$\CD
\TT@>w>>{}@>>>\TT\\
@Vf_{B'}VV    @.      @Vf_BVV\\
B'/U_{B'}@<\si<<(B\cap B')/(U_B\cap U_{B'})@>\si>>B/U_B  
\endCD$$
where the isomorphisms in the bottom row are induced by the obvious inclusions. We use 
this to identify $\WW$ with a subgroup of the group $\Aut(\TT)$ of automorphisms of 
the torus $\TT$. Let 
$$\WW^\bul=\{w\uD;w\in\WW,D\in G/G^0\}\sub\Aut(\TT).$$
This is a subgroup of $\Aut(\TT)$ normalizing $\WW$ since $\uD w=\e_D(w)\uD:\TT@>>>\TT$
for any $D\in G/G^0,w\in\WW$; here $\e_D$ is as in 26.2. 

Let $\la,\ra:\Hom(\kk^*,\TT)\T\Hom(\TT,\kk^*)@>>>\ZZ$ be the standard pairing. Define 
subsets $R,R^+$ of $\Hom(\TT,\kk^*)$ as follows. Let $B\in\cb$ and let $T$ be a maximal
torus of $B$. Consider the isomorphism $\TT@>\si>>T$ (composition of 
$f_B:\TT@>\si>>B/U_B$ with the obvious isomorphism $B/U_B@>\si>>T$). We require that 
the subset of $\Hom(T,\kk^*)$ corresponding to $R$ (resp. $R^+$) under this isomorphism
is the set of roots of $G^0$ with respect to $T$ (resp. the set of roots of $G^0$ with 
respect to $T$ such that the corresponding root subgroup is contained in $B$). Let 
$R^-=R-R^+$. For any $\a\in R$ there is a unique $\cha\in\Hom(\kk^*,\TT)$ and a unique 
$s_\a\in\WW$ such that $\la\cha,\a\ra=2$ and $s_\a(t)=t\cha(\a(t\i))$ for all 
$t\in\TT$. Then $s_\a^2=1$ and for $\cl\in\fs(\TT)$ we have 
$$\cl\cong s_\a^*\cl\ot\a^*(\cha^*\cl).\tag a$$
For $\cl\in\fs(\TT)$ let
$$R_\cl=\{\a\in R;\cha^*\cl\cong\bbq\}.$$
Pick $\k\in\Hom(\TT,\kk^*),\ce\in\fs(\kk^*)$ such that $\cl\cong\k^*\ce$, see 
28.1(b). Let $n\in\NN^*_\kk$ be the order of $\ce$ in $\ufs(\kk^*)$. We show:
$$R_\cl=\{\a\in R;\la\cha,\k\ra\in n\ZZ\}.$$ 
Indeed, for $\a\in R$ we have $\cha^*\cl=\cha^*\k^*\ce=(\k\circ\cha)^*\ce=f^*\ce$ where
$f:\kk^*@>>>\kk^*$ is $z\m z^{\la\cha,\k\ra}$. We now use the fact that, for $s\in\ZZ$,
the inverse image of $\ce$ under $\kk^*@>>>\kk^*,z@>>>z^s$ is $\bbq$ if and only if 
$s\in n\ZZ$.

Let $\WW^\bul_\cl=\{a\in\WW^\bul;(a\i)^*\cl\cong\cl\}$. Let $\WW_\cl$ be the subgroup 
of $\WW$ generated by $\{s_\a;\a\in R_\cl\}$. From (a) we see that
$$\WW_\cl\sub\WW^\bul_\cl.\tag b$$
Moreover, $\WW_\cl$ is a normal subgroup of $\WW^\bul_\cl$.

\subhead 28.4\endsubhead
In the remainder of this section we fix a connected component $D$ of $G$. Let 
$w\in\WW$. Let
$$\fZ_D^w=\{(B,\x);B\in\cb,\x\in U_B\bsl D/U_B,\po(B,gBg\i)=w\text{ for some/any }
g\in\x\}.$$

\subhead 28.5\endsubhead
Let $B^*\in\cb$ and let $T$ be a maximal torus of $B^*$. We set $U^*=U_{B^*}$,
$W_T=N_{G^0}T/T$. For any $w\in W_T$ we denote by $\dw$ a representative of $w$ in 
$N_{G^0}T$; we take $\dot 1=1$. We identify $W_T=\WW$ by 
$w\lra G^0\text{-orbit of }(B^*,\dw B^*\dw\i)$ and $T=\TT$ by 
$$t_1\lra t,t_1U^*=f_{B^*}(t).\tag a$$
Any $\a\in R$ becomes a root $\a$ of $G^0$ with respect to $T$ and $\cha$ becomes the 
corresponding coroot $\kk^*@>>>T$. We fix $d\in N_DB^*\cap N_DT$. For $w\in\WW$ we have
a diagram
$$T@<\ph<<\hfZ^w_D@>\r>>\fZ_D^w$$
where
$$\hfZ^w_D=\{(hU^*,g);hU^*\in G^0/U^*,g\in\dw dT\},$$
$$\r(hU^*,g)=(hB^*h\i, hU^*gU^*h\i),\ph(hU^*,g)=d\i\dw\i g.$$
Now $\ph$ is $T$-equivariant with respect to the $T$-action 

$t_0:(hU^*,g)\m(ht_0\i U^*,t_0gt_0\i)$
\nl
on $\hfZ^w_D$ and the $T$-action $t_0:\Ad(d\i\dw\i)(t_0)tt_0\i$ on $T$. Hence if 
$\cl_T\in\fs(T)$ satisfies $\Ad((\dw d)\i)^*\cl_T\cong\cl_T$, then $\ph^*\cl_T$ is a 
$T$-equivariant local system on $\hfZ^w_D$. (See 28.2(a).) Since $\r$ is a principal 
$T$-bundle, there is a well defined local system $\tcl_T$ on $\fZ^w_D$ such that 
$\r^*\tcl_T=\ph^*\cl_T$. 

\subhead 28.6\endsubhead
Let $w\in\WW$. Let $\cl\in\fs(\TT)$ be such that $w\uD\in\WW^\bul_\cl$. We associate to
$\cl$ a local system of rank $1$ $\tcl$ on $\fZ_D^w$ as follows. Let $B^*,U^*,T,d,\dw$ 
be as in 28.5. Using the identification $T=\TT$ in 28.5, we transport $\cl$ to a local 
system $\cl_T$ in $\fs(T)$. Then $\tcl_T$ is defined as in 28.5. Let $\tcl=\tcl_T$ (a 
local system on $\fZ_D^w$). We show:

(a) {\it the isomorphism class of $\tcl$ is independent of the choice of 
$B^*,T,\dw,d$.}
\nl
Let us replace $B^*,T,\dw,d$ by $xB^*x\i,xTx\i,x\dw x\i,xdx\i$ where $x\in G^0$. Define
${}'\hfZ^w_D,{}'\ph,{}'\r,{}'\cl_T,{}'\tcl_T$ in terms of this new choice in the same 
way as $\hfZ^w_D,\ph,\r,\cl_T,\tcl_T$ were defined in terms of $B^*,T,\dw,d$. We have a
commutative diagram 
$$\CD
\TT@>a>>T@<\ph<<\hfZ^w_D@>\r>>\fZ_D^w\\
@V=VV   @VbVV    @VcVV          @V=VV\\
\TT@>a'>>xTx\i@<{}'\ph<<{}'\hfZ^w_D@>{}'\r>>\fZ_D^w
\endCD$$
where $b(t)=xtx\i,c(hU^*,g)=(hx\i xU^*x\i,xgx\i)$, $a$ is given by 28.5(a) and $a'$ is 
the analogous isomorphism defined in terms of $xB^*x\i,xTx\i$ instead of $B^*,T$. Then 
${}'\cl_T=b^*\cl_T$, ${}'\ph^*{}'\cl_T=c^*\ph^*\cl_T$ and ${}'\tcl_T=\tcl_T$. Hence to
prove (a) it suffices to show that if $B^*,T,\dw,d$ are replaced by 
$B^*,T,\dw t_1,dt_2$ where $t_1,t_2\in T$ then the isomorphism class of $\tcl_T$ does 
not change. Note that $\hfZ^w_D,\r$ remain unchanged under the replacement above. 
However the map $\ph$ defined in terms of $B^*,T,\dw,d$ is replaced by the composition 
of $\ph$ with a left translation on $T$. It remains to use that the inverse image of 
$\cl_T$ under a left translation of $T$ is isomorphic to $\cl_T$, see 28.1(c).

\subhead 28.7\endsubhead
Let $J$ be a subset of $\II$. For any $B\in\cb$ we denote by $Q_{J,B}$ the unique 
parabolic in $\cp_J$ that contains $B$; we write $U_{J,B}$ instead of $U_{Q_{J,B}}$.
Let $w\in\WW$. Let $Z_{\em,J,D}^w$ be the set of all triples $(B,B',gU_{J,B})$ where 
$B\in\cb,B'\in\cb$, $g\in D,gBg\i=B'$, $\po(B,B')=w$. The map 
$$\z:Z_{\em,J,D}^w@>>>\fZ_D^w,\qua(B,B',gU_{J,B})\m(B,U_BgU_B)$$
is an affine space bundle.

Let $\cl\in\fs(\TT)$ be such that $w\uD\in\WW^\bul_\cl$. Then $\z^*\tcl$ is a local 
system on $Z_{\em,J,D}^w$ denoted again by $\tcl$.

\subhead 28.8\endsubhead
Let $\ww=(w_1,w_2,\do,w_r)$ be a sequence in $\WW$, let $[\ww]=w_1w_2\do w_r$ and let
$$\align Z^\ww_{\em,J,D}=&\{(B_0,B_1,\do,B_r,gU_{J,B_0});B_i\in\cb(i\in[0,r]),g\in D,
gB_0g\i=B_r,\\&\po(B_{i-1},B_i)=w_i(i\in[1,r])\}.\endalign$$ 
We define a morphism   
$$\z:Z^\ww_{\em,J,D}@>>>\fZ^{[\ww]}_D,(B_0,B_1,\do,B_r,gU_{J,B_0})\m
(B_0,U_{B_0}n_1n_2\do n_rnU_{B_0})$$
where $h_i\in G^0(i\in[1,r])$ are such that $B_i=h_iB_0h_i\i$, $h_0=1$, $T_0$ is a 
maximal torus of $B_0$, $n_i\in N_{G^0}T_0(i\in[1,r])$ are given by 
$h_{i-1}\i h_i\in U_{B_0}n_iU_{B_0}$ and $n\in N_DB_0\cap N_DT_0$ is given by 
$h_r\i g\in U_{B_0}n$. 

This is independent of the choices. (Another choice for $h_i,g,T_0$ must be of the form
$h'_i=h_iu_it_i,g'=gu',T'=uT_0u\i$ where $u_i\in U_{B_0}(i\in[1,r])$,
$t_i\in T_0(i\in[1,r])$, $u\in U_{B_0}$, $u'\in U_{J,B_0}$. Define $n'_i,n'$ in terms 
of this new choice in the same way as $n_i,n$ were defined in terms of the original 
choice. We have $n'_i=ut_{i-1}\i n_it_iu\i$ where $t_0=1$ and $n'=ut_r\i nu\i$. Hence 
$n'_1n'_2\do n'_rn'=un_1n_2\do n_rnu\i$ and
$$U_{B_0}n'_1n'_2\do n'_rn'U_{B_0}=U_{B_0}n_1n_2\do n_rnU_{B_0},$$
as required.) It is easy to see that $\z$ is an affine space bundle. Hence 
$Z^\ww_{\em,J,D}$ is smooth, connected.

Let $\cl\in\fs(\TT)$ be such that $[\ww]\uD\in\WW^\bul_\cl$. The inverse image under 
$\z$ of the local system $\tcl$ on $\fZ_D^{[\ww]}$ is a local system on 
$Z_{\em,J,D}^\ww$ denoted again by $\tcl$.

When $\ww$ has a single term $w$, we have $Z_{\em,J,D}^\ww=Z_{\em,J,D}^w$ and $\tcl$ 
defined above is the same as $\tcl$ defined in 28.7.

\subhead 28.9\endsubhead
Let $\ss=(s_1,s_2,\do,s_r)$ be a sequence in $\II\cup\{1\}$ and let
$$\align\bZ^\ss_{\em,J,D}=&\{(B_0,B_1,\do,B_r,gU_{J,B_0});B_i\in\cb(i\in[0,r]),g\in D,
gB_0g\i=B_r, \\&\po(B_{i-1},B_i)=1\text{ or }s_i(i\in[1,r])\}.\endalign$$
Let $\cj^0=\{j\in[1,r];s_j\in\II\}$. For any subset $\cj\sub\cj^0$ we consider the 
sequence $\ss_\cj=(s'_1,s'_2,\do,s'_r)$ in $\II\cup\{1\}$ given by $s'_i=s_i$ if 
$i\n\cj$ and $s'_i=1$ if $i\in\cj$; let $[\ss_\cj]=s'_1s'_2\do s'_r$. Then 
$Z^{\ss_\cj}_{\em,J,D}$ (see 28.8) is the locally closed subvariety of 
$\bZ^\ss_{\em,J,D}$ defined by the conditions 
$$B_{i-1}=B_i\text{ if }i\in\cj,\po(B_{i-1},B_i)=s_i\text{ if }i\n\cj.$$
The sets $Z^{\ss_\cj}_{\em,J,D} (\cj\sub\cj^0)$ form a partition of 
$\bZ^\ss_{\em,J,D}$. We have $\ss_\em=\ss$ and the corresponding piece 
$Z^{\ss_\em}_{\em,J,D}=Z^\ss_{\em,J,D}$ is open dense in $\bZ^\ss_{\em,J,D}$. Let 
$\cl\in\fs(\TT)$ be such that $[\ss]\uD\in\WW^\bul_\cl$. Let $\tcl$ be the local system
on $Z^\ss_{\em,J,D}$ defined as in 28.8. Let 
$\bcl=IC(\bZ^\ss_{\em,J,D},\tcl)\in\cd(\bZ^\ss_{\em,J,D})$. Let
$$\cj_\ss=\{j\in\cj^0;s_rs_{r-1}\do s_j\do s_{r-1}s_r\in\e_D(\WW_\cl)\}.$$

\proclaim{Lemma 28.10}$\bcl$ is a constructible sheaf on $\bZ^\ss_{\em,J,D}$. More 
precisely, $\bcl$ is a local system on the open subset 
$\cup_{\cj\sub\cj_\ss}Z^{\ss_\cj}_{\em,J,D}$ of $\bZ^\ss_{\em,J,D}$ and is $0$ on its 
complement; its restriction to $Z^{\ss_\cj}_{\em,J,D}$ (for $\cj\sub\cj_\ss$) is 
isomorphic to $\tcl$ (defined as in 28.8 in terms of $\ss_\cj$).
\endproclaim
For the statement above to make sense, we must verify that, if $\cj\sub\cj_\ss$, then 
$[\ss_\cj]\uD\in\WW^\bul_\cl$. This follows from $[\ss]\uD\in\WW^\bul_\cl$ and 
$\uD\i s_rs_{r-1}\do s_j\do s_{r-1}s_r\uD\in\WW^\bul_\cl$ for all $j\in\cj$. (See 
28.3(b).) 

Let $B^*,U^*,T,\ds_j,d$ be as in 28.5. Using the identification $T=\TT$ in 28.5 we
transport $\cl$ to a local system $\cl_T\in\fs(T)$. Let
$$\align\tZ^\ss=&\{(h_0,h_1,\do,h_r,g)\in G^0\T\do\T G^0\T D;\\&
h_{i-1}\i h_i\in B^*\ds_iB^*\cup B^*(i\in[1,r]),h_r\i gh_0\in N_GB^*\}.\endalign$$
The map $\tZ^\ss@>>>\bZ^\ss_{\em,J,D}$,
$$(h_0,h_1,\do,h_r,g)\m(h_0B^*h_0\i,h_1B^*h_1\i,\do,h_rB^*h_r\i,gU_{J,B_0})\tag a$$
is a locally trivial fibration with connected, smooth fibres. For $\cj\sub\cj^0$, the 
inverse image under (a) of the subvariety $Z^{\ss_\cj}_{\em,J,D}$ of 
$\bZ^\ss_{\em,J,D}$ is the subvariety $Z'{}^{\ss_\cj}$ of $\tZ^\ss$ defined by the 
conditions 
$$h_{i-1}\i h_i\in B^*\ds_i B^*(i\in[1,r]-\cj),h_{i-1}\i h_i\in B^*(i\in\cj).$$
It suffices to prove the statement analogous to that in the lemma for the inverse image
under (a) of $\bcl$. For $\cj\sub\cj^0$, define $\ps_\cj:Z'{}^{\ss_\cj}@>>>T$ by
$$(h_0,h_1,\do,h_r,g)\m d\i\dot{\ss}_\cj\i n_1n_2\do n_rn$$ where $n_i\in N_{G^0}T$ are
given by $h_{i-1}\i h_i\in U^*n_iU^*$ and $n\in N_GB^*\cap N_GT$ is given by 
$h_r\i gh_0\in U^*n$. (Here, we write $\dot{\ss}_\cj=\ds'_1\ds'_2\do\ds'_r$ for
$\ss_\cj=(s'_1,s'_2,\do,s'_r)$.) It suffices to prove the following statement:

{\it $IC(\tZ^\ss,\ps_\em^*\cl_T)$ is a local system on the open subset 
$\cup_{\cj\sub\cj_\ss}Z'{}^{\ss_\cj}$ of $\tZ^\ss$ and is $0$ on its complement; its
restriction to $Z'{}^{\ss_\cj}$ (for $\cj\sub\cj_\ss$) is isomorphic to 
$\ps_\cj^*(\cl_T)$.}
\nl
By the change of variables 
$$utd=h_r\i gh_0,y_i=h_{i-1}\i h_i(i\in[1,r-1]), y_r=h_{r-1}\i h_rt,$$
$\tZ^\ss$ becomes
$$\{(h_0,y_1,\do,y_r,u,t)\in G^0\T\do\T G^0\T U^*\T T;y_i\in B^*\ds_i B^*\cup B^*
(i\in[1,r])\}$$
and for $\cj\sub\cj^0$, $Z'{}^{\ss_\cj}$ becomes the subset of $\tZ^\ss$ defined by the
conditions
$$y_i\in B^*\ds_i B^*(i\in[1,r]-\cj),y_i\in B^*(i\in\cj).$$
Moreover, $\ps_\cj$ becomes 
$$(h_0,y_1,\do,y_r,u,t)\m d\i\dot{\ss}_\cj\i n_1n_2\do n_rd$$
where $n_i\in N_{G^0}T$ are given by $y_i\in U^*n_iU^*$. Since $h_0,u,t$ now play
passive roles, we can omit them. Thus, we set
$${}'\bZ^\ss=\{(y_1,\do,y_r)\in G^0\T\do\T G^0;y_i\in B^*\ds_i B^*\cup B^*(i\in[1,r])\}
$$
and, for $\cj\sub\cj^0$, we denote by ${}'Z^{\ss_\cj}$ the subset of ${}'\bZ^\ss$ 
defined by the conditions
$$y_i\in B^*\ds_i B^*(i\in[1,r]-\cj),y_i\in B^*(i\in\cj).$$
Define ${}'\ps_\cj:{}'Z^{\ss_\cj}@>>>T$ by 
$(y_1,\do,y_r)\m\dot{\ss}_\cj\i n_1n_2\do n_r$ where $n_i\in N_{G^0}T$ are given by 
$y_i\in U^*n_iU^*$. Let $\cl'=\Ad(d\i)^*\cl_T$. It suffices to prove the following 
statement:

{\it $IC({}'\bZ^\ss,{}'\ps_\em^*\cl')$ is a local system on the open subset 
$\cup_{\cj\sub\cj_\ss}{}'Z^{\ss_\cj}$ of ${}'\bZ^\ss$ and is $0$ on its complement; its
restriction to ${}'Z^{\ss_\cj}$ (for $\cj\sub\cj_\ss$) is isomorphic to 
$\ps_\cj^*\cl'$.}
\nl
The closures $\D_j$ of the subvarieties $\D_j^0={}'Z^{\ss_{\{j\}}} (j\in\cj^0)$ are 
smooth divisors with normal crossing in ${}'\bZ^\ss$. Using \cite{\CS, I,1.6} we see 
that it suffices to prove the following statement.

(b) {\it For $j\in\cj^0$, the monodromy of ${}'\ps_\em^*\cl'$ around the divisor 
$\D_j$ is trivial if and only if $j\in\cj_\ss$; if this condition is satisfied, then 
there exists a local system $\cf$ on the smooth variety $Z={}'Z^{\ss_\em}\cup\D_j^0$ 
such that $\cf|_{{}'Z^{\ss_\em}}\cong{}'\ps_\em^*\cl'$ and 
$\cf|_{\D_j^0}\cong{}'\ps_{\{j\}}^*\cl'$.}
\nl
Let $U_j$ be the root subgroup of $U^*$ with repect to $T$ such that 
$\ds_jU_j\ds_j\i\cap U^*=\{1\}$ and let $x:\kk@>\si>>U_j$ be an isomorphism. Let $\a$
be the root of $G^0$ with respect to $T$ such that $tx(a)t\i=x(\a(t)a)$ for all 
$t\in T$, $a\in\kk$. For $a\in\kk$ we set 
$y(a)=(\ds_1,\do,\ds_{j-1},\ds_jx(a)\ds_j\i,\ds_{j+1},\do,\ds_r)\in Z$. Then 
$y:\kk^*@>>>Z$ is a cross section to $\D_j^0$ in $Z$; we have $y(0)\in\D_j^0$, 
$y(a)\in Z-\D_j^0$ for $a\ne 0$. For $a\ne 0$ we have 
$\ds_jx(a)\ds_j\i\in U^*\ds_j\cha(a)U^*$. Hence
$${}'\ps_\em(y(a))=\dot{\ss}\i\ds_1\do\ds_{j-1}\ds_j\cha(a)\ds_{j+1}\do\ds_r
=\ds_r\i\do\ds_{j+1}\i\cha(a)\ds_{j+1}\do\ds_r=t_0\che\b_j(a)$$
where $\b_j$ is the root of $G^0$ with respect to $T$ that corresponds to the 
reflection $s_r\do s_{j+1}s_js_{j+1}\do s_r$ and $t_0$ is a fixed element of $T$. We 
see that $y^*{}'\ps_\em^*\cl'\cong\che\b_j^*\cl'$. Thus, 
$y^*{}'\ps_\em^*\cl'\cong\bbq$ if and only if $\che\b_j^*\cl'\cong\bbq$ that is, if
$s_r\do s_{j+1}s_js_{j+1}\do s_r\in\e_D(\WW_\cl)$ that is, if $j\in\cj_\ss$. This 
proves the first assertion of (b).

The second assertion of (b) involves only the component $G^0$ of $G$. Hence to prove 
it, we may assume that $G=G^0$. Let $\tG@>>>G$ be a surjective homomorphism of
connected reductive groups whose kernel is a central torus in $\tG$ and such that $\tG$
has simply connected derived group. The second assertion of (b) for $G$ follows from
the analogous assertion for $\tG$. Thus, we may assume that $G=G^0$ has simply
connected derived group. We may assume that $\cl'=\k^*\ce$ where 
$\k\in\Hom(T,\kk^*)$, $\ce\in\fs(\kk^*)$. Let $m$ be the order of $\ce$ in 
$\ufs(\kk^*)$. By assumption, we have $\che\b_j^*\cl'\cong\bbq$ hence 
$\la\che\b_j,\k\ra=mm_1$ with $m_1\in\ZZ$. Since $G$ has simply connected derived 
group, we can find $\k_1\in\Hom(T,\kk^*)$ such that $\la\che\b_j,\k_1\ra=m_1$. Then 
$\la\che\b_j,\k\k_1^{-m}\ra=0$. We have $(\k\k_1^{-m})^*\ce\cong\k^*\ce$. Hence 
replacing $\k$ by $\k\k_1^{-m}$, we may assume that $\la\che\b_j,\k\ra=0$. Then there 
is a unique homomorphism of algebraic groups $\c:B^*\ds_jB^*\cup B^*@>>>\kk^*$ such 
that 
$$\c(t)=\k(\ds_r\i\do\ds_{j+1}\i t\ds_{j+1}\do\ds_r)\text{ for all }t\in T.$$
We may assume that $\ds_j$ is in the derived subgroup of $B^*\ds_jB^*\cup B^*$. Then 
$\c(\ds_j)=1$. Define a morphism $\tf:Z@>>>\kk^*$ by
$$\align&\tf(y_1,\do,y_r)=\c(\ds_j\i\ds_{j-1}\i\do\ds_1\i n_1n_2\do n_{j-1}y_j
n_{j+1}\do n_r\ds_r\i\do\ds_{j+1}\i)\\&=\c(\ds_j\i\ds_{j-1}\i\do\ds_1\i n_1n_2\do 
n_{j-1}n_jn_{j+1}\do n_r\ds_r\i\do\ds_{j+1}\i),\endalign$$
where $n_i\in N_GT$ are given by $y_i\in U^*n_iU^*$. If 
$y_j\in B^*\ds_jB^*$, we have
$$\align&\tf(y_1,\do,y_r)=\k(\ds_r\i\do\ds_{j+1}\i\ds_j\i\ds_{j-1}\i\do\ds_1\i 
n_1n_2\do n_{j-1}n_jn_{j+1}\do n_r)\\&=\k({}'\ps_\em(y_1,\do,y_r)).\endalign$$
If $y_j\in B^*$, we have 
$$\align&\tf(y_1,\do,y_r)=\c(\ds_{j-1}\i\do\ds_1\i n_1n_2\do n_{j-1}n_j
n_{j+1}\do n_r\ds_r\i\do\ds_{j+1}\i)\\&=\k(\ds_r\i\do\ds_{j+1}\i\ds_{j-1}\i\do\ds_1\i 
n_1n_2\do n_{j-1}n_jn_{j+1}\do n_r)=\k({}'\ps_{\{j\}}(y_1,\do,y_r)).\endalign$$
Hence the local system $\cf=\tf^*(\ce)$ on $Z$ has the required properties and (b) is
verified. The lemma is proved.

\proclaim{Lemma 28.11}In the setup of 28.9 assume that $r\ge 2$, that $j\in[2,r]$ and 
$s_{j-1}=s_j\in\II$. Let $Z_1$ be the open subset of $Z^\ss_{\em,J,D}$ defined by 
$\po(B_{j-2},B_j)=s_j$. Let $\ss'=(s_1,s_2,\do,s_{j-1},s_{j+1},\do,s_r)$. Define 
$\d:Z_1@>>>Z^{\ss'}_{\em,J,D}$ by
$$(B_0,B_1,\do,B_r,gU_{J,B_0})\m(B_0,B_1,\do,B_{j-2},B_j,B_{j+1},\do,B_r,gU_{J,B_0}).$$
Let $\tcl$ be the local system on $Z^\ss_{\em,J,D}$ associated to $\cl$ as in 28.8; in 
the case where $j\in\cj_\ss$, let $\tcl'$ be the analogous local system on 
$Z^{\ss'}_{\em,J,D}$ associated to $\cl$. Let $\tcl_1$ be the restriction of $\tcl$ to 
$Z_1$. If $j\in\cj_\ss$, then $\tcl_1\cong\d^*\tcl'$. If $j\n\cj_\ss$, then 
$\d_!\tcl_1=0$.
\endproclaim
Consider the union $Z^{\ss_\em}_{\em,J,D}\cup Z^{\ss_{\{j\}}}_{\em,J,D}$ inside 
$\bZ^\ss_{\em,J,D}$. Recall that we may identify 
$Z^{\ss_\em}_{\em,J,D}=Z^\ss_{\em,J,D}$, 
$Z^{\ss_{\{j\}}}_{\em,J,D}=Z^{\ss'}_{\em,J,D}$. For
$$(B_0,B_1,\do,B_{j-2},B_j,B_{j+1},\do,B_r,gU_{J,B_0})\in Z^{\ss'}_{\em,J,D},$$
$$\align F=&\{(B_0,B_1,\do,B_{j-2},B,B_j,B_{j+1},\do,B_r,gU_{J,B_0});
\po(B_{j-2},B)=s_j,\\&\po(B,B_j)=1\text{ or }s_j\}\endalign$$
is a cross section to $Z^{\ss_{\{j\}}}_{\em,J,D}$ in $\bZ^\ss_{\em,J,D}$ which 
intersects $Z^{\ss_{\{j\}}}_{\em,J,D}$ in the point $\p$ defined by $B=B_j$. If 
$j\n\cj_\ss$, then the proof of Lemma 28.10 shows that the restriction of $\tcl$ to 
$F-\{\p\}\cong\kk^*$ is a local system in $\fs(\kk^*)$ not isomorphic to $\bbq$. Hence
$H^i_c(F-\{\p\},\tcl)=0$ for all $i$. Now $F-\{\p\}$ is the fibre of $\d$ at 
$(B_0,B_1,\do,B_{j-2},B_j,B_{j+1},\do,B_r,gU_{J,B_0})$. We see that the cohomology with
compact support of any fibre of $\d$ with coefficients in $\tcl$ is $0$. Thus, 
$\d_!\tcl=0$ as required. 

Assuming that $j\in\cj_\ss$, the same argument shows that the restriction of $\tcl_1$ 
to any fibre of $\d$ is $\bbq$. Hence $\tcl=\d^*\ce$ where $\ce$ is a local system of 
rank $1$ on $Z^{\ss'}_{\em,J,D}$. The proof of Lemma 28.10 shows that there exists a 
local system $\cf$ on $Z^\ss_{\em,J,D}\cup Z^{\ss'}_{\em,J,D}$ such that 
$\cf|_{Z^\ss_{\em,J,D}}=\tcl$ and $\cf|_{Z^{\ss'}_{\em,J,D}}=\tcl'$. Let $V$ be the 
open subset of $\bZ^\ss_{\em,J,D}$ defined by $\po(B_{i-1},B_i)=s_i$ if 
$i\in[1,r],i\ne j$ and $\po(B_{j-2},B_j)=s_j$. We have
$V=Z_1\cup Z^{\ss'}_{\em,J,D}$. Then $\cf^1:=\cf|_V$ is a local system on $V$ such that
$\cf^1|_{Z_1}=\tcl_1$ and $\cf^1|_{Z^{\ss'}_{\em,J,D}}=\tcl'$. Define 
$\tid:V@>>>Z^{\ss'}_{\em,J,D}$ by the same formula as $\d$. Then $\cf^1,\tid^*\ce$ are 
local systems on $V$ with the same restriction $\tcl_1$ on the open dense subset $Z_1$ 
of $V$. Hence $\cf^1\cong\tid^*\ce$. Since $\cf^1|_{Z^{\ss'}_{\em,J,D}}=\tcl'$, 
$\tid^*\ce|_{Z^{\ss'}_{\em,J,D}}=\ce$, we see that $\tcl'\cong\ce$. Thus, 
$\tcl\cong\d^*\tcl'$. The lemma is proved.

\subhead 28.12\endsubhead
Let $\e_D,Z_{J,D}$ be as in 26.2. Let $\cl\in\fs(\TT)$. For $w\in\WW$ we define 
$$\p:Z^w_{\em,J,D}@>>>Z_{J,D},(B,B',gU_{J,B})\m(Q_{J,B},Q_{\e_D(J),B'},gU_{J,B}).$$
If $w$ satisfies $w\uD\in\WW^\bul_\cl$, we set 
$K^{w,\cl}_{J,D}=\p_!\tcl\in\cd(Z_{J,D})$ where $\tcl$ is the local system on 
$Z^w_{\em,J,D}$ defined in 28.7.

For a sequence $\ww=(w_1,w_2,\do,w_r)$ in $\WW$ we define 
$\p_\ww:Z^\ww_{\em,J,D}@>>>Z_{J,D}$ by 
$$(B_0,B_1,\do,B_r,gU_{J,B_0})\m(Q_{J,B_0},Q_{\e_D(J),B_r},gU_{J,B_0}).\tag a$$
If $\ww$ satisfies $w_1w_2\do w_r\uD\in\WW^\bul_\cl$, we set
$K^{\ww,\cl}_{J,D}=\p_{\ww!}\tcl\in\cd(Z_{J,D})$ where $\tcl$ is the local system on 
$Z^\ww_{\em,J,D}$ defined in 28.8.

For a sequence $\ss=(s_1,s_2,\do,s_r)$ in $\II\cup\{1\}$ we define 
$\bpi_\ss:\bZ^\ss_{\em,J,D}@>>>Z_{J,D}$ by (a). If $\ss$ satisfies 
$s_1s_2\do s_r\uD\in\WW^\bul_\cl$, we set
$\bK^{\ss,\cl}_{J,D}=\bpi_{\ss!}\bcl\in\cd(Z_{J,D})$ where $\bcl$ is as in 28.9. Then

(b) $\bK^{\ss,\cl}_{J,D}\in\cd(Z_{J,D})$ {\it is a semisimple complex.}
\nl
This follows by applying the decomposition theorem \cite{\BBD, 5.4.5, 5.3.8} to the 
proper map $\bpi$.

\proclaim{Proposition 28.13} Let $\cl\in\fs(\TT)$ and let $A$ be a simple perverse 
sheaf on $Z_{J,D}$. The following conditions on $A$ are equivalent:

(i) $A\dsv K^{w,\cl}_{J,D}$ for some $w\in\WW$ such that $w\uD\in\WW^\bul_\cl$;

(ii) $A\dsv K^{\ww,\cl}_{J,D}$ for some $\ww=(w_1,w_2,\do,w_r)$ with $w_i\in\WW$, 
$w_1w_2\do w_r\uD\in\WW^\bul_\cl$;

(iii) $A\dsv K^{\ss,\cl}_{J,D}$ for some $\ss=(s_1,s_2,\do,s_r)$ with
$s_i\in\II\cup\{1\}$, $s_1s_2\do s_r\uD\in\WW^\bul_\cl$;

(iv) $A\dsv\bK^{\ss,\cl}_{J,D}$ for some $\ss=(s_1,s_2,\do,s_r)$ with 
$s_i\in\II\cup\{1\}$, $s_1s_2\do s_r\uD\in\WW^\bul_\cl$.

(v) $A\dsv\bK^{\ss,\cl}_{J,D}$ for some $\ss=(s_1,s_2,\do,s_r)$ with $s_i\in\II$, 
$s_1s_2\do s_r\uD\in\WW^\bul_\cl$.
\endproclaim
If in (ii), $\ww$ reduces to a single element $w$ then 
$K^{\ww,\cl}_{J,D}=K^{w,\cl}_{J,D}$. Thus, (i)$\imp$(ii). The implication 
(iii)$\imp$(ii) is trivial. We now prove the implication (ii)$\imp$(iii). Let 
$\ww=(w_1,w_2,\do,w_r)$ be a sequence in $\WW$ such that 
$w_1w_2\do w_r\uD\in\WW^\bul_\cl$ and let for some $i\in[1,r]$, $w'_i,w''_i$ be
elements of $\WW$ such that $w_i=w'_iw''_i$ and $l(w_i)=l(w'_i)+l(w''_i)$. Let 
$\ti\ww=(w_1,\do,w_{i-1},w'_i,w''_i,w_{i+1},\do,w_r)$. The map 
$$(B_0,B_1,\do,B_{r+1},gU_{J,B_0})\m
(B_0,B_1,\do,B_{i-1},B_{i+1},\do,B_{r+1},gU_{J,B_0})$$ 
defines an isomorphism $Z^{\ti\ww}_{J,D}@>\si>>Z^\ww_{J,D}$ compatible with the maps 
$\p_{\ti\ww},\p_\ww$ and with the local systems $\tcl$ defined on 
$Z^{\ti\ww}_{J,D},Z^\ww_{J,D}$ in terms of $\cl$ as in 28.8. Hence 
$$K^{\ti\ww,\cl}_{J,D}=K^{\ww,\cl}_{J,D}.\tag a$$
Applying (a) repeatedly we see that $K^{\ww,\cl}_{J,D}$ is equal to 
$K^{\ww',\cl}_{J,D}$ for some sequence $\ww'$ in $\II$. Thus, (ii)$\imp$(iii).

We prove the equivalence of (iii),(iv). Let $\ss=(s_1,s_2,\do,s_r)$ be a sequence in 
$\II\cup\{1\}$ such that $s_1s_2\do s_r\uD\in\WW^\bul_\cl$. Let $\bcl$ be as in 28.9. 

Define a sequence ${}^0Z\sps{}^1Z\sps\do$ of closed subsets of $\bZ^\ss_{\em,J,D}$ by 
$${}^iZ=\cup_{\cj\sub\cj^0;|\cj|\ge i}Z^{\ss_\cj}_{\em,J,D}.$$
Let 
$f^i:{}^iZ@>>>\bZ^\ww_{\em,J,D},f'{}^i:{}^iZ-{}^{i+1}Z@>>>\bZ^\ss_{\em,J,D}$ be the 
inclusions. The natural distinguished triangle in $\cd(Z_{J,D})$
$$(\bpi_\ss{}_!f'{}^i_!(f'{}^i)^*\bcl,\bpi_\ss{}_!f^i_!(f^i)^*\bcl,
\bpi_\ss{}_!f^{i+1}_!(f^{i+1})^*\bcl)$$
gives rise for any $i\ge 0$ to a long exact sequence in $\cm(Z_{J,D})$:
$$\align& \do@>>>{}^pH^{j-1}(\bpi_\ss{}_!f^{i+1}_!(f^{i+1})^*\bcl)
@>>>\op_{\cj\sub\cj_\ss;|\cj|=i}{}^pH^j(K^{\ss_\cj,\cl}_{J,D})\\&@>>>
{}^pH^j(\bpi_\ss{}_!f^i_!(f^i)^*\bcl)@>>>{}^pH^j(\bpi_\ss{}_!f^{i+1}_!(f^{i+1})^*\bcl)
@>>>\op\Sb\cj\sub\cj_\ss\\|\cj|=i\endSb{}^pH^{j+1}(K^{\ss_\cj,\cl}_{J,D})@>>>\do.\tag b
\endalign$$
Here we have used the isomorphism
$\bpi_\ss{}_!f'{}^i_!(f'{}^i)^*\bcl=\op_{\cj\sub\cj_\ss;|\cj|=i}K^{\ss_\cj,\cl}_{J,D}$
which follows from Lemma 28.10. Note that

($*$) $\bpi_\ss{}_!f^0_!(f^0)^*\bcl=\bK^{\ss,\cl}_{J,D},\bpi_\ss{}_!f^i_!(f^i)^*\bcl=0$
for $i$ large.
\nl
We set $m(\ss)=\sh(i\in[1,r];s_i\in\II)$. If $m(\ss)=0$ then 
$Z^\ss_{\em,J,D}=\bZ^\ss_{\em,J,D}$ and $K^{\ss,\cl}_{J,D}=\bK^{\ss,\cl}_{J,D}$. Hence 
in this case we have $A\dsv K^{\ss,\cl}_{J,D}$ if and only if 
$A\dsv\bK^{\ss,\cl}_{J,D}$. It suffices to verify the following statement.

(c) {\it Assume that $\ss$ satisfies $m(\ss)=m\ge 1$ and that for any sequence 
$\ss'=(s'_1,s'_2,\do,s'_r)$ in $\II\cup\{1\}$ with 
$s'_1s'_2\do s'_r\uD\in\WW^\bul_\cl$ and with $m(\ss')<m$ we have 
$A\not\dsv K^{\ss',\cl}_{J,D}$. Then $A\dsv K^{\ss,\cl}_{J,D}$ if and only if 
$A\dsv\bK^{\ss,\cl}_{J,D}$.}
\nl
Using (b) and our hypothesis we see that for any $i>0$ we have
$A\dsv\bpi_\ss{}_!f^i_!(f^i)^*\bcl$ if and only if 
$A\dsv\bpi_\ss{}_!f^{i+1}_!(f^{i+1})^*\bcl$. Applying this repeatedly for 
$i=N,N-1,\do,1$ (with $N$ large) we see that $A\not\dsv\bpi_\ss{}_!f^1_!(f^1)^*\bcl$. 
Using this, together with (b) we see that $A\dsv K^{\ss,\cl}_{J,D}$ if and only if 
$A\dsv\bpi_\ss{}_!f^0_!(f^0)^*\bcl$ that is, $A\dsv\bK^{\ss,\cl}_{J,D}$ (see $(*)$). 
This proves (c). The equivalence of (iii),(iv) is established.

The equivalence of (iv),(v) is obvious.

Let $\ss=(s_1,s_2,\do,s_r)$ be a sequence in $\II$ such that
$s_1s_2\do s_r\uD\in\WW^\bul_\cl$. Assume that $r\ge 2$ and that, for some $j\in[2,r]$ 
we have $s_{j-1}=s_j$. We have a partition $Z^\ss_{\em,J,D}=Z_1\cup Z_2$ where $Z_1$ 
(resp. $Z_2$) is the open (resp. closed) subset of $Z^\ss_{\em,J,D}$ defined by 
$\po(B_{j-2},B_j)=s_j$ (resp. by $B_{j-2}=B_j$). Let $\p_1,\p_2$ be the restrictions of
$\p_\ss$ to $Z_1,Z_2$. The natural distinguished triangle 
$(\p_{1!}\tcl,K^{\ss,\cl}_{J,D},\p_{2!}\tcl)$ in $\cd(Z_{J,D})$ (where the restrictions
of $\tcl$ from $Z^\ss_{\em,J,D}$ to $Z_1,Z_2$ are denoted again by $\tcl$) gives rise 
to a long exact sequence in $\cm(Z_{J,D})$:
$$\do@>>>{}^pH^i(\p_{1!}\tcl)@>>>{}^pH^i(K^{\ss,\cl}_{J,D})
@>>>{}^pH^i(\p_{2!}\tcl)@>>>{}^pH^{i+1}(\p_{1!}\tcl)@>>>\do.$$
Let $\ss'=(s_1,s_2,\do,s_{j-1},s_{j+1},\do,s_r)$,
$\ss''=(s_1,s_2,\do,s_{j-2},s_{j+1},\do,s_r)$. Then 
$$\d:(B_0,B_1,\do,B_r,gU_{J,B_0})\m(B_0,B_1,\do,B_{j-2},B_j,B_{j+1},\do,B_r,gU_{J,B_0})
$$
makes $Z_1$ into a locally trivial $\kk^*$-bundle over $Z^{\ss'}_{\em,J,D}$ and
$$(B_0,B_1,\do,B_r,gU_{J,B_0})\m(B_0,B_1,\do,B_{j-2},B_{j+1},\do,B_r,gU_{J,B_0})$$
makes $Z_2$ into a locally trivial affine line bundle over $Z^{\ss''}_{\em,J,D}$. The 
local system $\tcl$ on $Z_2$ is the inverse image of the local system $\tcl$ on 
$Z^{\ss''}_{\em,J,D}$ defined as in 28.8 in terms of $\cl$. By Lemma 28.11, if 
$j\in\cj_\ss$, the local system $\tcl$ on $Z_1$ is the inverse image under $\d$ of the 
local system $\tcl$ on $Z^{\ss'}_{\em,J,D}$ defined as in 28.8 in terms of $\cl$; if 
$j\n\cj_\ss$, then $\d_!\tcl=0$. It follows that
$$\p_{2!}\tcl=K^{\ss'',\cl}_{J,D}[[-1]]$$
and, if $j\n\cj_\ss$, we have $\p_{1!}\tcl=0$. If $j\in\cj_\ss$, we have a natural 
distinguished triangle $(\p_{1!}\tcl,K^{\ss',\cl}_{J,D}[[-1]],K^{\ss',\cl}_{J,D})$ in 
$\cd(Z_{J,D})$. Hence we have long exact sequences in $\cm(Z_{J,D})$:
$$\do@>>>{}^pH^i(\p_{1!}\tcl)@>>>{}^pH^i(K^{\ss,\cl}_{J,D})@>>>
{}^pH^{i-2}(K^{\ss'',\cl}_{J,D})(-1)@>>>{}^pH^{i+1}(\p_{1!}\tcl)@>>>\do,\tag d$$
$$\do@>>>{}^pH^i(\p_{1!}\tcl)@>>>{}^pH^{i-2}(K^{\ss',\cl}_{J,D})(-1)@>>>
{}^pH^i(K^{\ss',\cl}_{J,D})@>>>{}^pH^{i+1}(\p_{1!}\tcl)@>>>\do,\tag e$$
if $j\in\cj_\ss$, and isomorphisms 
$${}^pH^i(K^{\ss,\cl}_{J,D})@>\si>>{}^pH^{i-2}(K^{\ss'',\cl}_{J,D})(-1)\tag f$$
if $j\n\cj_\ss$.

We prove that (v)$\imp$(i). Assume that $A\dsv K^{\ss,\cl}_{J,D}$ where
$\ss=(s_1,s_2,\do,s_r)$ is as in (v). We may assume that $r$ is minimum possible. We 
want to show that (i) holds. Assume first that $l(s_1s_2\do s_r)<r$. We show that this 
contradicts the minimality of $r$. We can find $j\in[2,r]$ such that 
$l(s_js_{j+1}\do s_r)=r-j+1$ and $l(s_{j-1}s_j\do s_r)<r-j+2$. We can find 
$s'_j,s'_{j+1},\do,s'_r\in\II$ such that $s'_js'_{j+1}\do s'_r=s_js_{j+1}\do s_r=y$ and
$s'_j=s_{j-1}$. Let 
$$\uu'=(s_1,s_2,\do,s_{j-1},s'_j,s'_{j+1},\do,s'_r),\qua\uu''=(s_1,s_2,\do,s_{j-1},y).
$$
From (a) we have $K^{\ss,\cl}_{J,D}=K^{\uu',\cl}_{J,D}=K^{\uu'',\cl}_{J,D}$. Hence we 
may assume that $s_{j-1}=s_j$. If $j\n\cj_\ss$ then (f) shows that 
$A\dsv K^{\ss'',\cl}_{J,D}$; since the sequence $\ss''$ has $r-2$ terms, this 
contradicts the minimality of $r$. Assume now that $j\in\cj_\ss$. By the minimality of 
$r$ we have $A\not\dsv K^{\ss',\cl}_{J,D}$. From (e) it follows that 
$A\not\dsv\p_{1!}\tcl$. This, together with (d) shows that $A\dsv K^{\ss'',\cl}_{J,D}$.
This again contradicts the minimality of $r$. We see that $l(s_1s_2\do s_r)=r$. By (a),
we have $K^{\ss,\cl}_{J,D}=K^{w,\cl}_{J,D}$ where $w=s_1s_2\do s_r$ and the desired 
conclusion follows. Thus, we have (v)$\imp$(i). The proposition is proved.

\subhead 28.14\endsubhead
Let $A$ be a simple perverse sheaf on $Z_{J,D}$ and let $\cl\in\fs(\TT)$. We write 
$A\in\che Z_{J,D}^\cl$ if $A$ satisfies the equivalent conditions (i)-(v) in 28.13. We 
write $A\in\che Z_{J,D}$ if $A\in\che Z_{J,D}^\cl$ for some $\cl\in\fs(\TT)$; we then 
say that $A$ is a {\it parabolic character sheaf} on $Z_{J,D}$ (see \cite{\PCS}). 

In the case where $J=\II$ we identify $D=Z_{J,D}$ by $g\m(G^0,G^0,g)$; we write 
$A\in\che D^\cl,A\in\che D$ instead of $A\in\che Z_{J,D}^\cl,A\in\che Z_{J,D}$. We say 
that $A$ is a {\it character sheaf} on $D$ if $A\in\che D$. 

In the case where $J=\II$ we write $K^{w,\cl}_D,K^{\ss,\cl}_D,\bK^{\ss,\cl}_D$ instead 
of $K^{w,\cl}_{J,D}$, $K^{\ss,\cl}_{J,D}$, $\bK^{\ss,\cl}_{J,D}$.

\subhead 28.15\endsubhead
Let $A\in\che Z_{J,D}$. We can find $n\in\NN^*_\kk$ and $\cl\in\fs_n(\TT)$ such that 
$A\in\che Z_{J,D}^\cl$. We show that 

(a) {\it $A$ is equivariant for the action 
$$(z,x):(P,P',gU_P)\m(xPx\i,xP'x\i,xz^ngx\i U_{xPx\i})$$
of $H={}^D\cz_{G^0}^0\T G^0$ on $Z_{J,D}$.}
\nl
We can find $w\in\WW$ such that $w\uD\in\WW^\bul_\cl$ and $A\dsv\p_!\tcl$ where $\tcl$ 
is the local system on $Z^w_{\em,J,D}$ defined in 28.7 and 
$\p:Z^w_{\em,J,D}@>>>Z_{J,D}$ is as in 28.12. Now $H$ acts on $Z^w_{\em,J,D}$ and on 
$\fZ_D^w$ by 
$$(z,x):(B,B',gU_{J,B})\m(xBx\i,xB'x\i,xz^ngx\i U_{J,xBx\i}),$$
$$(z,x):(B,U_BgU_B)\m(xBx\i,U_{xBx\i}xz^ngx\i U_{xBx\i})$$ 
and $\p$ and $\z:Z_{\em,J,D}^w@>>>\fZ_D^w$ (see 28.7) are compatible with the 
$H$-actions. It then suffices to show that the local system $\tcl$ on $\fZ_D^w$ (see 
28.6) is $H$-equivariant. Let $B^*,U^*,T,d,\dw$ be as in 28.5. Let $\cl_T\in\fs(T)$ be
as in 28.6. We have $\cl_T\in\fs_n(T)$. It suffices to show that $\tcl_T$ (see 28.5) is
$H$-equivariant. Let $T@<\ph<<\hfZ^w_D@>\r>>\fZ_D^w$ be as in 28.5. Now $H$ acts on $T$
by $(z,x):t\m z^nt$ and on $\hfZ^w_D$ by $(z,x):(hU^*,g)\m(xhU^*,z^ng)$; note that 
$\ph,\r$ are compatible with the $H$-actions. Using the definitions we see that it 
suffices to show that $\cl_T$ is $H$-equivariant. This follows from the fact that 
$\cl_T\in\fs_n(T)$. This proves (a).

\subhead 28.16\endsubhead
Consider a sequence $\ss=(s_1,s_2,\do,s_r)$ in $\II\cup\{1\}$ with 
$s_1s_2\do s_r\uD\in\WW^\bul_\cl$. Let $\ss'=(s_2,s_3,\do,s_r,\e_D(s_1))$. We have
$s_2s_3\do s_r\e_D(s_1)\uD\in s_1\WW^\bul_\cl s_1=\WW^\bul_{\cl'}$ where 
$\cl'=s_1^*\cl\in\fs(\TT)$. We have an isomorphism 
$$Z^\ss_{\em,\II,D}@>\si>>Z^{\ss'}_{\em,\II,D},
(B_0,B_1,\do,B_r,g)\m(B_1,B_2,\do,B_r,gB_1g\i,g),$$
under which the local system $\tcl$ on $Z^\ss_{\em,\II,D}$ defined in 28.8 in terms of 
$\cl$ corresponds to the analogous local system $\tcl'$ on $Z^{\ss'}_{\em,\II,D}$ 
defined in terms of $\cl'$. It follows that

(a) $K^{\ss,\cl}_D=K^{\ss',\cl'}_D$.

\subhead 28.17\endsubhead
Let $\ss=(s_1,s_2,\do,s_r)$ be a sequence in $\II\cup\{1\}$ such that 
$s_1s_2\do s_r\uD\in\WW^\bul_\cl$. Let $m=\sh(i\in[1,r];s_i\ne 1)+\dim G$. For any 
$j\in\ZZ$ we have:
$${}^pH^j(\bK^{\ss,\cl}_D)\cong{}^pH^{2m-j}(\bK^{\ss,\cl}_D)\tag a$$ 
(in $\cm(D)$). This is a special case of the "relative hard Lefschetz theorem" 
\cite{\BBD, 6.2.10} applied to the projective morphism 
$\bpi_\ss:\bZ^\ss_{\em,\II,D}@>>>D$ and to the perverse sheaf $\bcl[m]$ on 
$\bZ^\ss_{\em,\II,D}$.

\head 29. Restriction functor for character sheaves\endhead
\subhead 29.1\endsubhead
Let $D$ be a connected component of $G$ and let $P$ be a parabolic of $G^0$ such that 
$N_DP\ne\em$. Let $L$ be a Levi of $P$. Let $G'=N_GP\cap N_GL$, a reductive group with 
$G'{}^0=L$. Let $D'=G'\cap D$, a connected component of $G'$. Let 
$\res_D^{D'}:\cd(D)@>>>\cd(D')$ be as in 23.3. Let $\a=\dim U_P$. We write 
$\e:\WW@>>>\WW$ instead of $\e_D:\WW@>>>\WW$ (see 26.2). 

In this section we begin the study of $\res_D^{D'}(A)$ where $A$ is a character sheaf 
on $D$. One the the main result of this section is that $\res_D^{D'}(A)$ is a direct 
sum of shifts of character sheaves on $D'$. (Here the words "shifts of" can be omitted 
but this will only come after further work in section 31.) The results in this section 
extend results in the connected case that appeared in \cite{\CS, I,\S3}. An obscure 
point in the proof in \cite{\CS, I, 3.5} (pointed out to me by J.G.M. Mars in 1985) is 
here eliminated following in part \cite{\LM}. 

Let $\cb^\da$ be the variety of Borel subgroups of $L$. We show that the canonical 
torus $\TT$ of $G^0$ (see 28.3) is canonically isomorphic to the analogously defined 
canonical torus $\TT^\da$ of $L$. We define a map 
$\prod_{B\in\cb}B/U_B@>>>\prod_{\b\in\cb^\da}\b/U_\b$ by 
$(g_BU_B)_{B\in\cb}\m(h_\b U_\b)$ where, for $\b\in\cb^\da$, $h_\b$ is the image of
$g_{\b U_P}$ under the obvious homomorphism $\b U_P@>>>\b$. This map restricts to a map
$\TT@>>>\TT^\da$ which is an isomorphism of tori. We use this isomorphism to identify 
$\TT^\da=\TT$. 

Similarly, the Weyl group $\WW$ of $G^0$  (see 26.1) contains the analogously defined
Weyl group $\WW^\da$ of $L$ as a subgroup. The imbedding $\WW^\da@>>>\WW$ is obtained
by associating to the $L$-orbit of $(\b,\b')\in\cb^\da\T\cb^\da$ the $G^0$-orbit of
$(\b U_P,\b' U_P)\in\cb\T\cb$. If $J$ is the subset of $\WW^\da$ analogous to the
subset $\II$ of $\WW$, then the imbedding $\WW^\da@>>>\WW$ restricts to an imbedding
$J\sub\II$. The length function of $\WW^\da$ is just the restriction of the length 
function of $\WW$. With the notation of 26.1, we have $\WW^\da=\WW_J$. Define 
$\po^\da:\cb^\da\T\cb^\da@>>>\WW_J$ in terms of $L$ in the same way as 
$\po:\cb\T\cb@>>>\WW$ was defined in terms of $G^0$.

For any Borel $B$ of $G^0$ we set $P^B=(P\cap B)U_P$, a Borel of $P$.

Let ${}^J\WW$ be as in 26.1. Then $y\m\WW_Jy$ is a bijection 
${}^J\WW@>\si>>\WW_J\bsl\WW$. We also have a bijection from the set of $P$-orbits on 
$\cb$ (for the conjugation action) to $\WW_J\bsl\WW$: the $P$-orbit of $B\in\cb$ 
corresponds to the $\WW_J$ coset of $\po(P^B,B)\in\WW$. Let $v(y)$ be the $P$-orbit on 
$\cb$ corresponding to $\WW_Jy,y\in\WW$.

If $y\in{}^J\WW$ and $s\in\II$, there are three possibilities for $ys$:

(i) $ys\in{}^J\WW$ and $l(ys)>l(y)$; then $v(y)\sub\ov{v(ys)}-v(ys)$.

(ii) $ys\in{}^J\WW$ and $l(ys)<l(y)$; then $v(ys)\sub\ov{v(y)}-v(y)$.

(iii) $ys\n{}^J\WW$; then $ysy\i\in J$ and $v(ys)=v(y)$.
\nl
For any $y\in\WW,g\in D$ we have $gv(y)g\i=v(\e(y))$.

Define a homomorphism $\p:N_GP@>>>G'$ by $\p(z\o)=z$ where $z\in G',\o\in U_P$ (see 
1.26).

\subhead 29.2\endsubhead 
Until the end of 29.9 we fix a sequence $\ss=(s_1,s_2,\do,s_r)$ in $\II$ and 
$\cl\in\fs(\TT)$ such that $s_1s_2\do s_r\uD\in\WW^\bul_\cl$. Let $\cj_\ss\sub[1,r]$ be
as in 28.9.

We shall write $\bZ^\ss$ instead of $\bZ^\ss_{\em,\II}$ (see 28.9). For any subset 
$\ct$ of $[1,r]$ let $Z^{\ss_\ct}$ be the subset of $\bZ^\ss$ defined by the conditions
$B_{i-1}=B_i$ if $i\in\ct$, $\po(B_{i-1},B_i)=s_i$ if $i\n\ct$. We have a partition 
$\bZ^\ss=\sqc_{\ct\sub[1,r]}Z^{\ss_\ct}$. Let
$$\bZ'=\{(B_0,B_1,\do,B_r,g)\in\bZ^\ss;g\in N_DP\}.$$
Define $\bpi':\bZ'@>>>D'$ by $(B_0,B_1,\do,B_r,g)\m\p(g)$.

Any sequence $\yy=(y_0,y_1,\do,y_r)$ in ${}^J\WW$ defines a locally closed subvariety 
$\bZ'_\yy$ of $\bZ'$:
$$\bZ'_\yy=\{(B_0,B_1,\do,B_r,g)\in\bZ^\ss;g\in N_DP,B_i\in v(y_i)(i\in[0,r])\}.$$
Clearly, $\bZ'_\yy=\em$ unless $\yy$ satisfies
$$y_i=y_{i-1} \text{ or }y_i=y_{i-1}s_i\text{ for all }i\in[1,r],\e(y_0)=y_r.\tag a$$
Let $\fii_\yy:\bZ'_\yy@>>>\bZ^\ss$ be the inclusion. Let $\bpi'_\yy:\bZ'_\yy@>>>D'$ be 
the restriction of $\bpi'$. 

\subhead 29.3\endsubhead
Until the end of 29.11 we fix a sequence $\yy$ satisfying 29.2(a). We set 
$$\align\dd(\yy)&=\a+\sh(i\in[1,r];y_{i-1}s_i\in{}^J\WW,l(y_{i-1}s_i)<l(y_{i-1})\\&
=\a+\sh(i\in[1,r];y_is_i\in{}^J\WW,l(y_is_i)<l(y_i)).\endalign$$
We show that these two definitions of $\dd(\yy)$ are equivalent. Let 

$c=\sh(i\in[1,r];y_{i-1}=y_i,y_is_i\in{}^J\WW,l(y_{i-1}s_i)<l(y_{i-1}))$,

$d=\sh(i\in[1,r];l(y_i)-l(y_{i-1})=-1)$,

$e=\sh(i\in[1,r];l(y_i)-l(y_{i-1})=1)$.
\nl
The two definitions of $\dd(\yy)$ are $\a+c+d$, $\a+c+e$. Hence it suffices to show
that $d=e$. Clearly, $l(y_r)-l(y_0)=e-d$. Since $y_r=\e(y_0)$, we have 
$l(y_r)=l(y_0)$ hence $e=d$, as required. An equivalent definition of $\dd(\yy)$ is
$$\align\dd(\yy)&=\a+\sh(i\in[1,r];v(y_is_i)\sub\ov{v(y_i)}-v(y_i))\\&=
\a+\sh(i\in[1,r];v(y_{i-1}s_i)\sub\ov{v(y_{i-1})}-v(y_{i-1})).\endalign$$
We define a sequence $(\ts_1,\ts_2,\do,\ts_r)$ in $\II\cup\{1\}$ by 
$$\ts_i=s_i\text{ if }y_{i-1}s_i\in\WW_Jy_i,\qua\ts_i=1\text{ if }y_{i-1}s_i\n\WW_Jy_i.
$$
Define $\tt=(t_1,t_2,\do,t_r)$ by $t_i=y_{i-1}\ts_iy_i\i$. Then $t_i\in J\cup\{1\}$ and
$$t_i=y_{i-1}s_iy_{i-1}\i\text{ if }y_{i-1}s_i\in\WW_Jy_{i-1},t_i=1\text{ if }
y_{i-1}s_i\n\WW_Jy_{i-1}.$$
Let 
$$\align\bZ^{\tt\da}&=\{(\b_0,\b_1,\do,\b_r,h)\in\cb^\da\T\do\T\cb^\da\T D';
\po^\da(\b_{i-1},\b_i)=1\text{ or }t_i,\\&\b_r=h\b_0h\i\}.\endalign$$
This is a variety like $\bZ^\ss_{\em,\II}$ in 28.9 with $G,D,\ss,\II$ replaced by
$G',D',\tt,J$. Define 
$$\r:\bZ'_\yy@>>>\bZ^{\tt\da},
(B_0,B_1,\do,B_r,g)\m(\p(P^{B_0}),\p(P^{B_1}),\do,\p(P^{B_r}),\p(g)).$$

\proclaim{Lemma 29.4} $\r$ is an iterated affine space bundle with fibres of dimension
$\dd(\yy)$.
\endproclaim
Let $F$ be the fibre of $\r$ over $(\b_0,\b_1,\do,\b_r,h)\in\bZ^{\tt\da}$. We show only
that 

(a) {\it $F$ is an iterated affine space bundle over a point and $\dim F=\dd(\yy)$}.
\nl
For any $k\in[0,r]$ let $F_k$ be the set of all sequences $(B_0,B_1,\do,B_k)$ in $\cb$ 
such that 
$$\po(B_{i-1},B_i)=1\text{ or }s_i(i\in[1,k]),B_i\in v(y_i)(i\in[0,k]), 
\p(P^{B_i})=\b_i(i\in[0,k]).$$
Let $F_{r+1}=F$. We have obvious maps 
$$F=F_{r+1}@>\x_{r+1}>>F_r@>\x_r>>F_{r-1}@>\x_{r-1}>>\do@>\x_1>>F_0.$$
It is easy to see that $F_0\cong\kk^{l(y_0)}$ and that
$\x_{r+1}:F_{r+1}@>>>F_r$ is an affine space bundle with fibres of dimension
$\a-l(y_r)=\a-l(y_0)$. Moreover, for $i\in[1,r]$, 

(b) $\x_i:F_i@>>>F_{i-1}$ {\it is an affine space bundle with fibres of dimension $1$ 
if $v(y_is_i)\sub\ov{v(y_i)}-v(y_i)$ and of dimension $0$, otherwise.}
\nl
Now (a) follows from (b). This completes the proof.

\subhead 29.5\endsubhead
Let $\cj=\{i\in[1,r];\ts_i=1\}$. We have $\cj=\cj_1\sqc\cj_2$ where 
$$\align&\cj_1=\{i\in\cj;v(y_{i-1})\sub\ov{v(y_{i-1}s_i)}-v(y_{i-1}s_i)\},\\&
\cj_2=\{i\in\cj;v(y_{i-1}s_i)\sub\ov{v(y_{i-1})}-v(y_{i-1})\}.\endalign$$
Let $\ck^0=\{i\in[1,r];t_i\in J\}$. We have 
$$\cj\cap\ck^0=\em.$$
Indeed, if $i\in\cj\cap\ck^0$ then $\ts_i=1$ hence 
$v(y_i)\ne v(y_{i-1}s_i),v(y_i)=v(y_{i-1})$, $v(y_{i-1})\ne v(y_{i-1}s_i)$ and 
$t_i\in J,\ts_i=1$ hence $v(y_{i-1}s_i)=v(y_{i-1})$, a contradiction.

We show:

(a) {\it If $i\in\ck^0$ then $v(y_is_i)=v(y_i)$. If $i\in\cj_1$ then 
$v(y_i)\sub\ov{v(y_is_i)}-v(y_is_i)$. If $i\in\cj_2$ then 
$v(y_is_i)\sub\ov{v(y_i)}-v(y_i)$.}
\nl
Assume first that $i\in\ck^0$. Then $t_i\ne 1$, $v(y_{i-1}s_i)=v(y_{i-1})$, $s_i\ne 1$.
If $v(y_i)=v(y_{i-1})$ we get $v(y_is_i)=v(y_i)$; if $v(y_i)=v(y_{i-1}s_i)$ we get 
again $v(y_is_i)=v(y_i)$.

Assume next that $i\in\cj_1$. Then $v(y_{i-1})\sub\ov{v(y_{i-1}s_i)}-v(y_{i-1}s_i)$, 
$\ts_i=1$, $v(y_i)\ne v(y_{i-1}s_i)$ hence $v(y_i)=v(y_{i-1})$ and 
$v(y_i)\sub\ov{v(y_is_i)}-v(y_is_i)$.

Finally, assume that $i\in\cj_2$. Then $v(y_{i-1}s_i)\sub\ov{v(y_{i-1})}-v(y_{i-1})$, 
$\ts_i=1$, $v(y_i)\ne v(y_{i-1}s_i)$ hence $v(y_i)=v(y_{i-1})$ and 
$v(y_is_i)\sub\ov{v(y_i)}-v(y_i)$.

\subhead 29.6\endsubhead
For any subset $\ck\sub\ck^0$ let 
$$\align Z^{\tt_\ck\da}&=\{(\b_0,\b_1,\do,\b_r,h)\in\bZ^{\tt\da};
\po^\da(\b_{i-1},\b_i)=t_i(i\in[1,r]-\ck),\\&\b_{i-1}=\b_i(i\in\ck)\}.\endalign$$
We shall write $Z^{\tt_\da}$ instead of $Z^{\tt_\em\da}$. We have 
$\bZ^{\tt\da}=\sqc_{\ck\sub\ck^0}Z^{\tt_\ck\da}$. Hence
$\bZ'_\yy=\sqc_{\ck\sub\ck^0}Z'_{\yy,\ck}$ where $Z'_{\yy,\ck}=\r\i(Z^{\tt_\ck\da})$.
We show that for $(B_0,B_1,\do,B_r,h)\in\bZ'_{\yy,\ck}$, conditions (i),(ii) below are 
equivalent:

(i) $\ck\cup\cj_1\sub\{i\in[1,r];B_{i-1}=B_i\}\sub\ck\cup\cj$;

(ii) $\{i\in\ck^0;\p(P^{B_{i-1}})=\p(P^{B_i})\}=\ck$.
\nl
Assume that (i) holds. If $i\in\ck$, then by (i) we have $B_{i-1}=B_i$ hence 
$\p(P^{B_{i-1}})=\p(P^{B_i})$. Conversely, let $i\in\ck^0$ be such that 
$\p(P^{B_{i-1}})=\p(P^{B_i})$. Using 29.4(b) we see that $B_{i-1}=B_i$ (since 
$v(y_is_i)\not\sub\ov{v(y_i)}-v(y_i)$, by 29.5(a)). Using (i) we see that 
$i\in\ck\cup\cj$. Since $\ck^0\cap\cj=\em$ we deduce that $i\in\ck$. We see that 
(i)$\imp$(ii).

Assume that (ii) holds. If $i\in\ck$ then, by (ii), we have 
$\p(P^{B_{i-1}})=\p(P^{B_i})$; using (b) we see that $B_{i-1}=B_i$ (since
$v(y_is_i)\not\sub\ov{v(y_i)}-v(y_i)$, by (c)). If $i\in\cj_1$ then $\ts_i=1$,
$v(y_i)\ne v(y_{i-1}s_i)$ hence $v(y_i)=v(y_{i-1})$ and $v(y_{i-1})\ne v(y_{i-1}s_i)$. 
Then $t_i=1$ hence $\p(P^{B_{i-1}})=\p(P^{B_i})$. Using (b) we see that $B_{i-1}=B_i$ 
(since $v(y_is_i)\not\sub\ov{v(y_i)}-v(y_i)$, by (c)). Thus, the first inclusion in (i)
holds. Conversely, if $i\in\ck^0,B_{i-1}=B_i$ then $\p(P^{B_{i-1}})=\p(P^{B_i})$ and 
using (ii) we see that $i\in\ck$. If $i\in[1,r]-\ck^0,B_{i-1}=B_i$, then $t_i=1$, 
$v(y_{i-1}s_i)\ne v(y_{i-1})$, $v(y_i)=v(y_{i-1})$ hence $v(y_i)\ne v(y_{i-1}s_i)$, 
$\ts_i=1$ hence $i\in\cj$. Thus the second inclusion in (i) holds. We see that 
(ii)$\imp$(i).

The equivalence of (i),(ii) can be also formulated as follows:
$$Z'_{\yy,\ck}=\sqc_{\cj';\cj'\sub\cj_2}Z^{\ss_{\ck\cup\cj_1\cup\cj'}}\cap\bZ'_\yy.
\tag a$$

\subhead 29.7\endsubhead
Let $(\b_0,\b_1,\do,\b_r,h)\in\bZ^{\tt_\ck\da}$ and let $F,F_k,\x_k$ be as in 29.4.
From 29.6(a) we see that 
$$F=\sqc_{\cj';\cj'\sub\cj_2}(F\cap Z^{\ss_{\ck\cup\cj_1\cup\cj'}}).$$
For $i\in[1,r]$, let $F^i=\{(B_0,B_1,\do,B_r,g)\in F;B_{i-1}=B_i\}$. We show:

(a) {\it For $i\in\cj_2$, $F^i$ is a smooth hypersurface in $F$. For $i\in\ck\cup\cj_1$
we have $F^i=F$. For $i\in[1,r]-(\ck\cup\cj)$ we have $F^i=\em$.}
\nl
If $F^i\ne\em$ then, using $F\sub\cup_{\cj'\sub\cj_2}Z^{\ss_{\ck\cup\cj_1\cup\cj'}}$,
we see that $i\in\ck\cup\cj_1\cup\cj'$ for some $\cj'\sub\cj_2$; thus, 
$i\in\ck\cup\cj$. In the rest of the proof we assume that $i\in\ck\cup\cj$.

For each $k\in[i,r]$ let $F_k^i$ be the set of all $(B_0,B_1,\do,B_k)\in F_k$ such that
$B_{i-1}=B_i$. Let $F_{r+1}^i=F^i$. From the definitions we see that for 
$k\in[i+1,r+1]$ we have a cartesian diagram
$$\CD F^i_k@>>>F_k\\
     @VVV   @V\x_k VV\\
     F^i_{k-1}@>>>F_{k-1}\endCD$$
where the horizontal maps are inclusions. 

Assume first that $i\in\cj_2$. Using the cartesian diagram above, it suffices to show 
that $F^i_i$ is a smooth hypersurface in $F_i$. From 29.5(a) we see that 
$v(y_is_i)\sub\ov{v(y_i)}-v(y_i)$; hence $\x_i:F_i@>>>F_{i-1}$ is an affine line bundle
(see 29.4(b)). It suffices to show that $\x_i$ restricts to an isomorphism 
$F_i^i@>\si>>F_{i-1}$. Let $(B_0,B_1,\do,B_{i-1})\in F_{i-1}$. It suffices to show that
$(B_0,B_1,\do,B_{i-1},B_{i-1})\in F_i$. Hence it suffices to show that 
$v(y_{i-1})=v(y_i)$ and $\b_{i-1}=\b_i$. Since $i\in\cj$ we have $\ts_i=1$ hence 
$v(y_i)\ne v(y_{i-1}s_i)$ hence $v(y_i)=v(y_{i-1})$ and $t_i=1$. Since 
$\po^\da(\b_{i-1},\b_i)=1\text{ or }t_i$, we see that $\b_{i-1}=\b_i$, as required. 

Assume next that that $i\in\ck\cup\cj_1$. Using the cartesian diagram above, it 
suffices to show that $F^i_i=F_i$. From 29.5(a) we see that 
$v(y_is_i)\not\sub\ov{v(y_i)}-v(y_i)$; hence $\x_i:F_i@>>>F_{i-1}$ is an isomorphism 
(see 29.4(b)). It suffices to show that $\x_i$ restricts to an isomorphism 
$F_i^i@>\si>>F_{i-1}$. If $i\in\cj_1$ this is shown exactly as in the first part of the
proof. Assume now that $i\in\ck$. We have $t_i\ne 1$ hence 
$v(y_{i-1}s_i)=v(y_{i-1})=v(y_i)$. From the definitions we have $\b_{i-1}=\b_i$. Hence 
$F_i^i@>\si>>F_{i-1}$ as in the first part of the proof. This proves (a).

\proclaim{Lemma 29.8}The map $\r_1:\bZ'_\yy\cap Z^{\ss_\cj}@>>>Z^{\tt\da}$ (restriction
of $\r$) is an iterated affine space bundle.
\endproclaim
Let $(\b_0,\b_1,\do,\b_r,h)\in Z^{\tt\da}$. We show only that the fibre $\bF$ of $\r_1$
at \lb $(\b_0,\b_1,\do,\b_r,h)$ is an iterated affine space bundle over a point and
$$\dim\bF=\a+\sh(i\in[1,r]-\cj;v(y_is_i)\sub\ov{v(y_i)}-v(y_i)).$$
For any $k\in[0,r]$ let 
$$\align\bF_k=&\{(B_0,B_1,\do,B_k)\in\cb^{k+1};\po(B_{i-1},B_i)
=s_i(i\in[1,k],i\n\cj),\\&B_{i-1}=B_i(i\in[1,k]\cap\cj),B_i\in v(y_i) 
(i\in[0,k]),\p(P^{B_i})=\b_i(i\in[0,k])\}.\endalign$$ 
Let $\bF_{r+1}=\bF$. We have obvious maps 
$$\bF=\bF_{r+1}@>>>\bF_r@>>>\bF_{r-1}@>>>\do@>>>\bF_0.$$
It is easy to see that $\bF_0\cong\kk^{l(y_0)}$ and that $\bF_{r+1}@>>>\bF_r$ is an 
affine space bundle with fibres of dimension $\a-l(y_r)=\a-l(y_0)$. Moreover, for 
$i\in[1,r]$, $\bF_i@>>>\bF_{i-1}$ is an affine space bundle with fibres of dimension 
$1$ if $v(y_is_i)\sub\ov{v(y_i)}-v(y_i)$, $i\n\cj$ and of dimension $0$, otherwise. 
This completes the proof.

\subhead 29.9\endsubhead
For $k\in[1,r]$ we set 
$$\us_{r,k}=s_rs_{r-1}\do s_k\do s_{r-1}s_r,\qua
\tius_{r,k}=\ts_r\ts_{r-1}\do\ts_{k+1}s_k\ts_{k+1}\do\ts_{r-1}\ts_r.$$
Let $\G$ be a subgroup of $W$ such that 

$k\in[1,r],\ts_k=1\imp\us_{r,k}\in\G$.
\nl
We show:

(a) {\it For $i\in[1,r]$ we have $\us_{r,i}\in\G$ if and only if $\tius_{r,i}\in\G$.}
\nl
We argue by induction on $r-i$. If $r-i=0$ we have $\us_{r,i}=\tius_{r,i}$ so the 
result is obvious. Assume now that $r-i\ge 1$. We have 

$k\in[1,r-1],\ts_k=1\imp\us_{r-1,k}\in s_r\G s_r$.
\nl
By the induction hypothesis we have 
$\us_{r-1,i}\in s_r\G s_r$ if and only if $\tius_{r-1,i}\in s_r\G s_r$. If $\ts_r=s_r$
then $\us_{r-1,i}=s_r\us_{r,i}s_r$, $\tius_{r-1,i}=s_r\tius_{r,i}s_r$. Hence we have 
$\us_{r,i}\in\G$ if and only if $\tius_{r,i}\in\G$. If $\ts_r=1$ then 
$\us_{r-1,i}=s_r\us_{r,i}s_r$, $\tius_{r-1,i}=\tius_{r,i}$. Hence we have 
$\us_{r,i}\in\G$ if and only if $\tius_{r,i}\in\G$. (We use that $s_r\in\G$.) This 
proves (a).

We show:

(b) {\it If $i\in[1,r]$ and $t_i\ne 1$ then}
$t_rt_{r-1}\do t_i\do t_{r-1}t_r=y_r\tius_{r,i}y_r\i$.
\nl
We argue by induction on $r-i$. If $r-i=0$ we have $t_r\ne 1$ hence
$t_r=y_{r-1}\ts_ry_r\i=y_{r-1}s_ry_{r-1}\i$. We see that 
$y_{r-1}=y_r\ts_rs_r=y_rs_r\ts_r$ and $t_r=y_rs_ry_r\i$ as required.

Assume now that $r-i\ge 1$. By the induction hypothesis the left hand side of the
equality in (b) is $t_ry_{r-1}\tius_{r-1,i}y_{r-1}\i t_r$ and the right hand side is
$y_r\ts_r\tius_{r-1,i}\ts_ry_r\i$. It then suffices to show that 
$t_ry_{r-1}=y_r\ts_r$; this folows from the definitions since $t_r=t_r\i$. This proves
(b).

\subhead 29.10\endsubhead
We set $y=y_0$. In the case where $\cj\sub\cj_\ss$ we set 
${}^y\cl=\Ad(y\i)^*\cl\in\fs(\TT)$ and we show:
$$t_1t_2\do t_r\uD\in\WW^\bul_{{}^y\cl}.\tag a$$
For $i\in[1,r]$ define $u_i\in\WW$ by $u_i=\e\i(s_rs_{r-1}\do s_i\do s_{r-1}s_r)$ if
$i\in\cj$ and by $u_i=1$ if $i\n\cj$. If $i\in\cj$ we have $i\in\cj_\ss$ hence
$u_i\in\WW^\bul_\cl$. Then 
$$\ts_1\ts_2\do\ts_r\uD=s_1s_2\do s_r\uD u_1u_2\do u_r\in\WW^\bul_\cl.$$
We have $\ts_1\ts_2\do\ts_r=y\i t_1t_2\do t_ry_r$ hence 
$y\i t_1t_2\do t_ry_r\uD\in\WW^\bul_\cl$. Since $y_r=\e(y)$ we have 
$y_r\uD=\uD y$ hence $y\i t_1t_2\do t_r\uD y\in\WW^\bul_\cl$ and (a) follows.

Using (a) we can define a constructible sheaf ${}^y\bcl$ on $\bZ^{\tt\da}$ and a 
complex $\bK^{\tt,{}^y\cl}_{D'}\in\cd(D')$ in terms of $\tt,{}^y\cl,G'$ in the same way
as $\bcl$ on $\bZ^\ss$ and $\bK^{\ss,\cl}_D\in\cd(D)$ are defined in 28.12 in terms of 
$\ss,\cl,G$. 

\proclaim{Lemma 29.11}(a) If $\cj\not\sub\cj_\ss$ then $\bpi'_{\yy!}\fii_\yy^*\bcl=0$.

(b) If $\cj\sub\cj_\ss$ then 
$\bpi'_{\yy!}\fii_\yy^*\bcl=\bK^{\tt,{}^y\cl}_{D'}[[-\dd(\yy)]]$.
\endproclaim
Let $\bpi_\tt:\bZ^{\tt\da}@>>>D'$ be the obvious projection. We have 
$\bpi'_{\yy!}=\bpi_{\tt!}\r_!$ and $\r_!\r^*({}^y\bcl)={}^y\bcl[[-\dd(\yy)]]$ (we use
29.4). Hence it suffices to prove:

(a${}'$) {\it If $\cj\not\sub\cj_\ss$ then $\r_!(\fii_\yy^*\bcl)=0$.}

(b${}'$) {\it If $\cj\sub\cj_\ss$ then $\fii_\yy^*\bcl\cong\r^*({}^y\bcl)$.}
\nl
We prove (a${}'$). Let $F$ be the fibre of $\r$ over a point of 
$Z^{\tt_\ck\da}(\ck\sub\ck^0)$. We must show that $H^*_c(F,\bcl|_F)=0$. (We write 
$\bcl$ instead of $\fii_\yy^*\bcl$.) If $\ck\cup\cj_2\not\sub\cj_\ss$ then $\bcl|_F=0$ 
(we use 29.6(a) and Lemma 28.10) and the desired vanishing follows. Assume now that 
$\ck\cup\cj_2\sub\cj_\ss$ but $\cj\not\sub\cj_\ss$. Using 29.6(a) and Lemma 28.10 we 
see that $\bcl|_F$ is a local system on 
$$\sqc_{\cj';\cj'\sub\cj_2;\cj'\sub\cj_\ss}(F\cap Z^{\ss_{\ck\cup\cj_1\cup\cj'}})$$
and is zero elsewhere. Hence $\bcl|_F$ is a local system on 
$F-\cup_{j\in\cj_2-\cj_\ss}F^j$ ($F^j$ as in 29.7) and is zero on 
$\cup_{j\in\cj_2-\cj_\ss}F^j$. Let $i$ be the largest number in $\cj_2-\cj_\ss$. It 
suffices to show that for any $(B_0,B_1,\do,B_{i-1})\in F_{i-1}$ (see 29.4) we have
$H^*_c(F',\bcl|_{F'})=0$ where $F'$ is the fibre of the obvious map $F@>>>F_{i-1}$ at 
$(B_0,B_1,\do,B_{i-1})$. Let $F'{}^i=F'\cap F^i$. If $B_{j-1}=B_j$ for some 
$j<i,j\in\cj_2-\cj_\ss$ then $F'\sub F^j$ and $\bcl|_{F'}=0$; the desired vanishing
follows. Thus we may assume that $\po(B_{j-1},B_j)=s_j$ for all 
$j<i,j\in\cj_2-\cj_\ss$. Then $\bcl|_{F'-F'{}^i}$ is a local system and
$\bcl|_{F'{}^i}=0$. Let $F''$ be the fibre of $\x_i:F_i@>>>F_{i-1}$ (see 29.4) at 
$(B_0,B_1,\do,B_{i-1})$. Let $F''{}^i=\{(B_0,B_1,\do,B_{i-1},B_{i-1})\}$, a point on
the affine line $F''$. Let $u:F'@>>>F''$ be the restriction of the obvious map
$F@>>>F_i$. Then $u$ is an iterated affine space bundle (see 29.4(a)), 
$F'{}^i=u\i(F''{}^i)$, and there is a well defined local system $\ce$ on 
$F''-F''{}^i\cong\kk^*$ such that $\ce\in\fs(\kk^*)$, $\bcl|_{F'-F'{}^i}=u^*(\ce)$. 
Then $H^*_c(F'-F'{}^i,\bcl)\cong H^*_c(\kk^*,\ce)$ and it suffices to show that 
$\ce\not\cong\bbq$. It also suffices to show that $\bcl|_{F'-F'{}^i}$ has non-trivial 
monodromy around the smooth hypersurface $F'{}^i$ of $F'$. This is the same as the 
monodromy of $\bcl|_{Z^\ss}$ around the hypersurface $Z^{\ss_{\{i\}}}$. This monodromy 
is non-trivial by 28.10(b). This proves (a${}'$).

We prove (b${}'$). We define $\WW_{J,{}^y\cl}$ in terms of $G',\WW_J,{}^y\cl$ in the 
same way as $\WW_\cl$ was defined in terms of $G^0,\WW,\cl$. Let 
$$\align\cj_\tt
&=\{i\in[1,r];t_i\in J,t_rt_{r-1}\do t_i\do t_{r-1}t_r\in\e(\WW_{J,{}^y\cl})\}\\&
=\{i\in[1,r];t_i\in J,t_rt_{r-1}\do t_i\do t_{r-1}t_r\in\e(\WW_{{}^y\cl})\}.\endalign$$
(The two definitions coincide since $t_rt_{r-1}\do t_i\do t_{r-1}t_r$ is a reflection 
in $\WW_J$.) We show:
$$\cj_\tt=\cj_\ss\cap\ck^0.\tag c$$ 
Using 29.9(b) it suffices to show that for $i\in[1,r]$ such that $t_i\ne 1$, we have
$$y_r\ts_r\do\ts_{i+1}s_i\ts_{i+1}\do\ts_ry_r\i\in
\e(\WW_{{}^y\cl})\lra s_rs_{r-1}\do s_i\do s_{r-1}s_r\in\e(\WW_\cl).$$
Using 29.9(a), we see that it suffices to show that 
$\e\i(y_r\i)\WW_{{}^y\cl}\e\i(y_r)=\WW_\cl$ or that ${}^y\cl=\Ad(\e\i(y_r\i))^*\cl$. 
This follows from the definitions using $y_r=\e(y)$.

Using Lemma 28.10 for $G'$ instead of $G$, we see that ${}^y\bcl$ is a local system on 
$Z^1=\cup_{\ck\sub\cj_\tt}Z^{\tt_\ck\da}$ and is zero on its complement in
$\bZ^{\tt\da}$. Using Lemma 28.10 and 29.6(a) we see that $\fii_\yy^*\bcl$ is a local 
system on
$$Z^2=\cup_{\cj'\sub\cj_2;\ck\sub\cj_\ss}Z^{\ss_{\ck\cup\cj_1\cup\cj'}}\cap\bZ'_\yy
=\cup_{\ck\sub\cj_\tt}Z'_{\yy,\ck}$$
and is zero on its complement in $\bZ'_\yy$. (We have used (c).) Since $Z^2$ is an 
iterated affine space bundle over $Z^1$ (via $\r$) and the restriction of 
$\fii^*_\yy\bcl$ to any fibre of $\r:Z^2@>>>Z^1$ is a local system of rank $1$ with 
finite monodromy of order invertible in $\kk$ (hence it is $\bbq$) we see that 
$\fii^*_\yy\bcl|_{Z^2}=\r^*\ce$ for a well defined local system $\ce$ of rank $1$ on
$Z^1$. It suffices to show that $\ce\cong{}^y\bcl|_{Z^1}$. Since $Z^1$ is smooth and
$Z^{\tt\da}$ is open dense in $Z^1$, it suffices to show that 
$\ce|_{Z^{\tt\da}}\cong{}^y\bcl|_{Z^{\tt\da}}$. Let $Z^3=Z^{\ss_\cj}\cap\bZ'_\yy$. This
is a closed subsetof the open subset $\r\i(Z^{\tt\da})$ of $Z_2$. Since the restriction
of $\r$ is an iterated affine space bundle $\r_1:Z^3@>>>Z^{\tt\da}$ (see 29.8), it 
suffices to show that $\r_1^*(\ce|_{Z^{\tt\da}})\cong\r_1^*({}^y\bcl|_{Z^{\tt\da}})$. 
Since $\r_1^*(\ce|_{Z^{\tt\da}})=\bcl|_{Z^3}$, it suffices to show that 
$\bcl|_{Z^3}\cong\r_1^*({}^y\bcl|_{Z^{\tt\da}})$. Using Lemma 28.10, once for $G$ and
once for $G'$, we see that $\bcl|_{Z^3}=\tcl|_{Z^3}$, 
${}^y\bcl|_{Z^{\tt\da}}={}^y\tcl$ where $\tcl$ (on $Z^{\ss_\cj}$) is defined as in 28.8
in terms of $G,\ss_\cj,\cl$ and ${}^y\tcl$ (on $Z^{\tt\da}$) is defined analogously in 
terms of $G',\tt,{}^y\cl$. Thus it suffices to show that 
$$\tcl|_{Z^3}\cong\r_1^*({}^y\tcl).\tag d$$
To prove (d), we choose $B^*,T,d$ and $\dw$ (for $w\in\WW$) as in 28.5, in such a way 
that $B^*\sub P,T\sub L$. We have necessarily $d\in D'$. Let 
$\b^\da=\p(B^*)\in\cb^\da$. Let $U^*=U_{B^*},U^\da=U_{\b^\da}$. Let 
$$\align\cz=&\{(h_0U^*,h_1U^*,\do,h_rU^*,g)\in(G^0/U^*)^{r+1}\T D;
h_{i-1}\i h_i\in B^*\dts_iB^* \\&\text{ for }i\in[1,r],h_r\i gh_0\in N_GB^*,
g\in N_GP,h_i\in P\dy_iU^*\text{ for }i\in[0,r]\},\endalign$$
$$\align\cz'&=\{(h'_0U^\da,h'_1U^\da,\do,h'_rU^\da,g')\in(L/U^\da)^{r+1}\T D';\\&
h'_{i-1}{}\i h'_i\in\b^\da\dt_i\b^\da\text{ for }i\in[1,r],
h'_r{}\i g'h'_0\in N_G\b^\da\}.\endalign$$
Define $\z:\cz@>>>\cz'$ by 
$$(h_0U^*,h_1U^*,\do,h_rU^*,g)\m(h'_0U^\da,h'_1U^\da,\do,h'_rU^\da,g')$$
where $h_i\in p_i\dy_iU^*,p_i\in P,h'_i=\p(p_i),g'=\p(g)$. 
(We show that $h'_iU^\da$ is well defined or equivalently that $p_iU^*$ is well
defined. It suffices to show that $p\dy_iU^*=p'\dy_iU^*,p,p'\in P\imp p'U^*=pU^*$. It 
also suffices to show that $P\cap\dy_i U^*\dy_i\i\sub U^*$. Since $y_i\in{}^J\WW$ and 
$B^*\sub P$ we have $P^{\dy_iB^*\dy_i\i}=B^*$. Hence $P\cap\dy_iB^*\dy_i\i\sub B^*$. 
Thus $P\cap\dy_i U^*\dy_i\i$ is contained in the set of unipotent elements of $B^*$ 
that is, in $U^*$.) We have a commutative diagram
$$\CD
\cz@>\z>>\cz'\\
@VaVV    @Va'VV\\
Z^3@>\r_1>>Z^{\tt\da}
\endCD$$
where $a:(h_0U^*,h_1U^*,\do,h_rU^*,g)\m(h_0B^*h_0\i,h_1B^*h_1\i,\do,h_rB^*h_r\i,g)$,
$$a':(h'_0U^\da,h'_1U^\da,\do,h'_rU^\da,g')\m
(h'_0\b^\da h'_0{}\i,h'_1\b^\da h'_1{}\i,\do,h'_r\b^\da h'_r{}\i,g').$$
Since $a$ is a locally trivial fibration with smooth connected fibres, to prove (d) it 
suffices to prove that $a^*(\tcl|_{Z^3})\cong a^*\r_1^*({}^y\tcl)$ or that
$a^*(\tcl|_{Z^3})\cong\z^*a'{}^*({}^y\tcl)$. Define $\x:\cz@>>>T$ by
$$(h_0U^*,h_1U^*,\do,h_rU^*,g)\m d\i(\ds'_1\ds'_2\do\ds'_r)\i n_1n_2\do n_rn$$
where $n_i\in N_{G^0}T$ are given by $h_{i-1}\i h_i\in U^*n_iU^*$ and 
$n\in N_GB^*\cap N_GT$ is given by $h_r\i gh_0\in U^*n$. Define $\x_1:\cz'@>>>T$ by
$$(h'_0U^\da,h'_1U^\da,\do,h'_rU^\da,g')\m 
d\i(\dt_1\dt_2\do\dt_r)\i\bn_1\bn_2\do\bn_r\bn$$ 
where $\bn_i\in N_LT$ are given by $h'_{i-1}{}\i h'_i\in U^\da\bn_iU^\da$ and 
$\bn\in N_G\b^\da\cap N_GT$ is given by $h'_r{}\i g'h'_0\in U^\da\bn$. From the 
definitions we have
$$a^*(\tcl|_{Z^3})=\x^*\cl_T,a'{}^*({}^y\tcl)=\x_1^*\Ad(\dy\i)^*\cl_T$$
where $\cl_T\in\fs(T)$ is as in 28.6. Therefore it suffices to show that
$$\x^*\cl_T\cong\z^*\x_1^*\Ad(\dy\i)^*\cl_T.$$
Define $\x':\cz@>>>T$ by
$$(h_0U^*,h_1U^*,\do,h_rU^*,g)\m d\i(\dt_1\dt_2\do\dt_r)\i\bn_1\bn_2\do\bn_r\bn$$
where $h_i\in p_i\dy_iU^*$, $p_i\in P$, $p_{i-1}\i p_i\in U^*\bn_iU^*$,
$p_r\i gp_0\in U^*\bn$, $\bn_i\in N_LT$, $\bn\in N_GB^*\cap N_GT$. Then $\x'=\x_1\z$ 
and it suffices to show that 
$$\x^*\cl_T\cong\x'{}^*\Ad(\dy\i)^*\cl_T.$$
Using 28.1(c) it suffices to show that there exists $t\in T$ such that
$$t\x(z)=\Ad(\dy\i)\x'(z)\text{ for all }z\in\cz.$$
Let $(h_0,h_1,\do,h_r,g)\in(G^0)^{r+1}\T D$ be such that
$z=(h_0U^*,h_1U^*,\do,h_rU^*,g)\in\cz$. We define $p_i,n_i,n,\bn_i,\bn$ in terms of 
$h_i$ as in the definition of $\x,\x'$. From $h_r\i gh_0\in U^*n$, 
$p_r\i gp_0\in U^*\bn$, we deduce $\dy_r\i p_r\i gp_0\dy\in U^*n$, 
$U^*\bn\dy U^*=U^*\dy_r nU^*$ hence $\bn\dy=\dy_rn$. We show that 
$$\bn_i=\dy_{i-1}n_i\dy_i\i\text{ for any }i\in[1,r].$$
From $h_{i-1}\i h_i\in U^*n_iU^*$, $p_{i-1}\i p_i\in U^*\bn_iU^*$, we deduce
$\dy_{i-1}\i p_{i-1}\i p_i\dy_i\in U^*n_iU^*$ hence
$\dy_{i-1}un_iu'\dy_i\i\in U^*\bn_iU^*$ for some $u,u'\in U^*$. Assume first that 
$t_i\ne 1$. Then 
$$y_i=y_{i-1},y_{i-1}s_i=t_iy_i,l(y_{i-1}s_i)=l(t_iy_i)=l(y_{i-1})+1=l(y_i)+1,$$
hence 
$$\dy_{i-1}un_iu'\in U^*\dy_{i-1}n_iU^*, U^*\bn_iU^*\dy_i\sub U^*\bn_i\dy_iU^*.$$
Thus, $U^*\dy_{i-1}n_iU^*=U^*\bn_i\dy_iU^*$ and $\dy_{i-1}n_i=\bn_i\dy_i$, as required.
Next, assume that $t_i=1,\ts_i\ne 1$. Then $y_i=y_{i-1}s_i\ne y_{i-1}$, $\bn_i\in T$.
If $l(y_{i-1}s_i)=l(y_{i-1})+1$, then $\dy_{i-1}un_iu'\in U^*\dy_{i-1}n_iU^*$ and 
$U^*\bn_iU^*\dy_i\sub U^*\bn_i\dy_iU^*$ so that $U^*\dy_{i-1}n_iU^*=U^*\bn_i\dy_iU^*$ 
and $\dy_{i-1}n_i=\bn_i\dy_i$, as required. If $l(y_{i-1}s_i)=l(y_{i-1})-1$, then 
$l(s_iy_i\i)=l(y_i\i)+1$. We have $un_iu'\dy_i\i\in U^*n_i\dy_i\i U^*$ and 
$\dy_{i-1}\i U^*\bn_iU^*\sub U^*\dy_{i-1}\i\bn_iU^*$ so that
$U^*n_i\dy_i\i U^*=U^*\dy_{i-1}\i\bn_iU^*$ and $n_i\dy_i\i=\dy_{i-1}\i\bn_i$, as 
required. Finally, assume that $\ts_i=1$. Then $t_i=1,y_{i-1}s_i\ne y_{i-1}=y_i$, 
$n_i\in T$, $\bn_i\in T$. We have $\dy_{i-1}un_iu'\in U^*\dy_{i-1}n_iU^*$ and 
$U^*\bn_iU^*\dy_i\sub U^*\bn_i\dy_iU^*$ so that $U^*\dy_{i-1}n_iU^*=U^*\bn_i\dy_iU^*$ 
and $\dy_{i-1}n_i=\bn_i\dy_i$, as required. 

We have
$$\align&\Ad(\dy\i)\x'(z)=\dy\i d\i(\dt_1\dt_2\do\dt_r)\i\bn_1\bn_2\do\bn_r\bn\dy
\\&=\dy\i d\i(\dt_1\dt_2\do\dt_r)\i(\dy n_1\dy_1\i)(\dy_1n_2\dy_2\i)\do
(\dy_{r-1}n_r\dy_r\i)\dy_r n\\&=\dy\i d\i(\dt_1\dt_2\do\dt_r)\i\dy n_1n_2\do n_rn\\&
=td\i(\dts_1\dts_2\do\dts_r)\i n_1n_2\do n_rn=t\x(z)\endalign$$
where
$$t=\dy\i d\i(\dt_1\dt_2\do\dt_r)\i\dy(\dts_1\dts_2\do\dts_r)d.$$
We have $t\in T$. (Equivalently, $y\ts_1\ts_2\do\ts_r=t_1\do t_ry_r$, which is clear 
from the definitions.) This completes the proof of  (d) hence that of (b${}'$). The 
lemma is proved.

\subhead 29.12\endsubhead
We consider the sequence $Z_0\sub Z_1\sub\do$ of closed subsets of $\bZ'$ defined by
$Z_i=\cup_{\yy;c(\yy)\le i}\bZ'_\yy$ where $\yy$ is a sequence $(y_0,y_1,\do,y_r)$ of 
elements in ${}^J\WW$ satisfying 29.2(a) and $c(\yy)=\sum_{i\in[0,r]}\dim v(y_i)$. Let 
$k_i:Z_i@>>>\bZ'(i\ge 0)$ and $k'_i:Z_i-Z_{i-1}@>>>\bZ'(i\ge 1)$ be the inclusions. For
any $i\ge 1$, the natural distinguished triangle
$$(\bpi'_!k'_{i!}k'_i{}^*\bcl,\bpi'_!k_{i!}k_i^*\bcl,\bpi'_!(k_{i-1})_!k_{i-1}^*\bcl)$$
in $\cd(D')$ gives rise to a long exact sequence in $\cm(D')$:
$$\align&\do@>>>{}^pH^{j-1}(\bpi'_!(k_{i-1})_!k_{i-1}^*\bcl)@>\d>>
\op_{\yy;c(\yy)=i}{}^pH^j(\bpi'_{\yy!}\fii_\yy^*\bcl)@>>>\\&
{}^pH^j(\bpi'_!k_{i!}k_i^*\bcl)@>>>{}^pH^j(\bpi'_!(k_{i-1})_!k_{i-1}^*\bcl)@>\d>>\do.
\tag a\endalign$$
We now prove the following result.

\proclaim{Lemma 29.13}(a) The maps $\d$ in 29.12(a) are zero.

(b) For $i\ge 0$, $\bpi'_!k_{i!}k_i^*\bcl\in\cd(D')$ is a semisimple complex; it is 
isomorphic to $\op_{\yy;c(\yy)\le i}\bpi'_{\yy!}\fii_\yy^*\bcl$.

(c) $\bpi'_!\bcl\in\cd(D')$ is a semisimple complex; it is isomorphic to
$\op_\yy\bpi'_{\yy!}\fii_\yy^*\bcl$.
\endproclaim
(c) is a special case of (b), for large $i$. Assuming that (a) and the first assertion 
of (b) are proved, we prove the second assertion of (b) as follows. Since both 
complexes in question are semisimple (see 29.11 and 28.12(b)), it suffices to show that
they have the same ${}^pH^j$ for any $j$. Using (a) we see that 29.12(a) decomposes 
into short exact sequences of semisimple objects in $\cm(D')$. Hence
$${}^pH^j(\bpi'_!k_{i!}k_i^*\bcl)\cong{}^pH^j(\bpi'_!(k_{i-1})_!k_{i-1}^*\bcl)\op
\op_{\yy;c(\yy)=i}{}^pH^j(\bpi'_{\yy!}\fii_\yy^*\bcl).$$
This proves the desired equality for ${}^pH^j$ by induction on $i$. (The case where
$i=0$ is trivial.) 

It remains to prove (a) and the first assertion of (b). By general principles, we may 
assume that $\kk$ is an algebraic closure of a finite field $\FF_q$, that $G,P,D$ are 
defined over $\FF_q$ and that $G^0$ is split over $\FF_q$. By taking $\FF_q$ large
enough, we may assume that 29.12(a) and the isomorphisms in 29.11(a),(b) are realized 
in the category of mixed perverse sheaves with $\cl$ pure of weight $0$. Now 
$\bK^{\tt,{}^{y_0}\cl}_{D'}$ in 29.11(b) is pure of weight $0$ (by Deligne's theorem 
\cite{\DE, 6.2.6}) since it is a direct image under a proper map of ${}^{y_0}\bcl$ 
which is pure of weight $0$; after applying to it $[[-\dd(\yy)]]$, it remains pure of 
weight $0$, see \cite{\BBD, 6.1.4}. Hence by 29.11, 
$\bpi'_\yy{}_!(\fii_\yy^*\bcl)$ is pure of weight $0$; it follows that

(d) $\op_{\yy;c(\yy)=i}{}^pH^j(\bpi'_{\yy!}\fii_\yy^*\bcl)$ is pure of weight $j$.
\nl
We now show by induction on $i$ that ${}^pH^j(\bpi'_!k_{i!}k_i^*\bcl)$ is pure of
weight $j$ for any $i$. For $i=0$ this follows from (d). If we assume that this holds 
for $i-1$ where $i\ge 1$ then the statement for $i$ follows from the statement for 
$i-1$ and 29.12(a) (using (d)); we also use the following fact: if $K_1@>>>K_2@>>>K_3$ 
is an exact sequence of mixed perverse sheaves with $K_1,K_3$ pure of weight $j$, then
$K_2$ is pure of weight $j$. Using \cite{\BBD, 5.4.4} it follows that 
$\bpi'_!k_{i!}k_i^*\bcl$ is pure of weight $0$. Using the decomposition theorem
\cite{\BBD, 5.4.5, 5.3.8} it follows that $\bpi'_!k_{i!}k_i^*\bcl$ is a semisimple
complex. The vanishing of $\d$ in 29.12(a) follows from the fact that $\d$ is a 
morphism between two pure perverse sheaves of different weights. The lemma is proved.

\proclaim{Proposition 29.14}In $\cd(D')$ we have
$\res_D^{D'}(\bK^{\ss,\cl}_D)\cong\op_\yy\bK^{\tt,{}^{y_0}\cl}_{D'}[[-\dd(\yy)]]$ where
$\yy$ runs over all sequences satisfying 29.2(a) such that 
$\{i\in[1,r],\ts_i=1\}\sub\cj_\ss$ and $\ts_i,\tt$ are defined in terms of $\yy$ as in 
29.3. In particular, $\res_D^{D'}(\bK^{\ss,\cl}_D)$ is a direct sum of shifts of 
character sheaves on $D'$.
\endproclaim
From the definitions we have $\res_D^{D'}(\bK^{\ss,\cl}_D)=\bpi'_!\bcl(\a)$. The result
follows from 29.13(c) and 29.11.

\proclaim{Proposition 29.15} Let $A$ be a character sheaf on $D$. Then 
$\res_D^{D'}A\in\cd(D')$ is a direct sum of shifts of character sheaves on $D'$.
\endproclaim
We can find $\ss,\cl$ as in 29.2 such that $A\dsv\bK^{\ss,\cl}_D$. Using 28.12(b) we 
see that for some $j\in\ZZ$, $A[-j]$ is a direct summand of $\bK^{\ss,\cl}_D$. Hence 
$\res_D^{D'}A[-j]$ is a direct summand of $\res_D^{D'}(\bK^{\ss,\cl}_D)$, which is a
semisimple complex by 29.14. It follows that $\res_D^{D'}A$ is a semisimple complex. 
Now ${}^pH^i(\res_D^{D'}A)$ is a direct summand of 
${}^pH^{i+j}(\res_D^{D'}(\bK^{\ss,\cl}_D))$ which, by 29.14, is a direct sum of 
character sheaves on $D'$. Hence ${}^pH^i(\res_D^{D'}A)$ is a direct sum of character 
sheaves on $D'$. This completes the proof.

\head 30. Admissibility of character sheaves\endhead
\subhead 30.1\endsubhead
In this section we fix a connected component $D$ of $G$. We write $\e:\WW@>>>\WW$ 
instead of $\e_D:\WW@>>>\WW$ (see 26.2). 

\proclaim{Lemma 30.2} Let $H={}^D\cz_{G^0}^0\T G^0$. Let $A$ be a simple perverse sheaf
on $D$ which is cuspidal (see 23.3). Assume that there exists $n\in\NN^*_\kk$ such that
$A$ is equivariant for the $H$-action $(z,x):g\m xz^ngx\i$ on $D$. Let $Z=\supp A$, 
$m=\dim Z$. There exists a unique pair $(S,\ce)$ where $S$ is an isolated stratum $S$ 
of $D$ and $\ce$ is an irreducible cuspidal local system $\ce\in\cs(S)$ (up to 
isomorphism) such that $A[-m]=IC(\bS,\ce)$ extended by $0$ on $D-\bS$.
\endproclaim
The intersections of $Z$ with the various strata of $D$ form a finite partition of $Z$ 
into locally closed subsets. Since $Z$ is irreducible, one of these intersections is
open dense in $Z$. Thus there exists $(L,S)\in\AA$ such that $S\sub D$ and 
$Y_{L,S}\cap Z$ is open dense in $Z$. Let $P$ be a parabolic of $G^0$ with Levi $L$ 
such that $S\sub N_GP$. Let $a=\dim U_P$. We can find an open dense smooth subset $V$ 
of $Z$ and an irreducible local system $\ce$ on $V$ such that $A=IC(Z,\ce)[m]$ extended
by $0$ on $D-Z$. Replacing if necessary $V,\ce$ by 
$V\cap Y_{L,S},\ce|_{V\cap Y_{L,S}}$, we may assume that $V\sub Y_{L,S}$. For any 
$h\in H$, the $h$-translate ${}^hV$ of $V$ is an open dense smooth subset of $Z$. Hence
$V'=\cup_h{}^hV$ is an open dense smooth subset of $Z$. Since $V\sub Y_{L,S}$ and 
$Y_{L,S}$ is $H$-stable, we have ${}^hV\sub Y_{L,S}$ for $h\in H$ 
hence $V'\sub Y_{L,S}$. Now $A'=A[-m]|_V$ is an $H$-equivariant intersection cohomology
complex on $V'$ such that $A'|_V$ is a local system and $A'|_{{}^hV}$ is automatically 
a local system for any $h\in H$. Since $\cup_h{}^hV$ is an open covering of $V'$, we
see that $A'$ is a local system on $V'$. Replacing $V,\ce$ by $V',A'$, we see that we 
may assume in addition that $V$ is $H$-stable and $\ce$ is an $H$-equivariant local 
system on $V$. Define $f:G^0\T(V\cap S^*)@>>>V$ by $(y,g)\m ygy\i$. Then $f$ is 
surjective since $V\sub Y_{L,S}$. Moreover, $f$ is a principal bundle with group
$\G=\{x\in N_{G^0}L;xSx\i=S\}$ which acts on $G^0\T(V\cap S^*)$ by 
$x:(y,g)\m(yx\i,xgx\i)$. (We show this only at the level of sets. It suffices to show
that, if $(y,g),(y',g')$ are elements of $V\cap S^*$ such that $ygy\i=y'g'y'{}\i$, then
the element $x=y'{}\i y\in G^0$ satisfies $xLx\i=L,xSx\i=S$. We have $xgx\i=g'$. Since 
$g\in S^*,g'\in S^*$, we have $L=L(g)=L(g')$, see 3.9, and $L(g')=xL(g)x\i$ hence
$xLx\i=L$. Since $xSx\i,S$ are strata of $N_GL$ with a common element $g$, we must have
$xSx\i=S$, as required.) Since $V$ is irreducible and $\G$ is of pure dimension 
$\dim L$, it follows that $V\cap S^*$ is non-empty, of pure dimension $m-2a$: we have
$$\dim(V\cap S^*)+\dim G^0=\dim V+\dim\G=\dim V+\dim L=m+\dim G^0-2a.$$
Let $g\in V\cap S^*$. Let $U'$ be the orbit of $g$ under $U_P$-conjugation. Since $U'$ 
is an orbit of an action of a unipotent group on an affine variety, $U'$ is closed in 
$D$. We have $U'\sub gU_P$. (Indeed if $x\in U_P$ then $xgx\i=g(g\i xg)x\i$ and 
$g\i xg\in U_P$ since $g\in N_GP$.) The isotropy group $U_{P,g}$ of $g$ in $U_P$ is 
contained in 
$$U_P\cap Z_G(g)\sub U_P\cap Z_G(g_s)\sub U_P\cap Z_G(g_s)^0$$
(the last inclusion follows from 1.11). Since $g\in S^*$, we have $Z_G(g_s)^0\sub L$ 
hence $U_{P,g}\sub U_P\cap L=\{1\}$. Thus, $U_{P,g}=\{1\}$. We see that 
$\dim U'=\dim U_P$. Since $U'$ is closed in $gU_P$, we have $U'=gU_P$. Since $V$ is 
stable under $U_P$-conjugation and $U'$ is the $U_P$-orbit of $g\in V$, it follows that
$U'\sub V$. Thus, $gU_P\sub V$. Now $\ce|_{gU_P}$ is a $U_P$-equivariant local system 
(for the conjugation action of $U_P$ which has trivial isotropy group). It follows that
$\ce|_{gU_P}\cong\bbq^c$ for some $c\ge 1$. Hence $H^{2a}_c(gU_P,\ce)\ne 0$.
Equivalently, $H^{2a-m}_c(gU_P,A)\ne 0$. 

For any $i\in\ZZ$, we denote by $\cx^i$ the set of all $U_P$-cosets $R$ in $N_DP$ such 
that $H^i_c(R,A)\ne 0$. Then, for any $g\in V\cap S^*$, we have $gU_P\in\cx^{2a-m}$.
The map
 $$V\cap S^*@>>>N_DP/U_P,g\m gU_P$$
is injective: if $g,g'\in V\cap S^*$ and $gU_P=g'U_P$ then 
$$g\i g'\in(N_GP\cap N_GL)\cap U_P=\{1\},$$ 
see 1.26, hence $g=g'$. We see that $\dim\cx^{2a-m}\ge\dim(V\cap S^*)$ hence 
$\dim\cx^{2a-m}\ge m-2a$. Thus,
$$\dim(\supp\ch^{2a-m}(\res_D^{D'}A))\ge m-2a$$
where $D'=N_DP\cap N_DL$. If $P\ne G^0$ then our assumption that $A$ is cuspidal gives 
$$\dim(\supp\ch^{2a-m}(\res_D^{D'}A))<m-2a,$$
a contradiction. Thus, $P=G^0,L=G^0$ and $S$ must be an isolated stratum of $D$, so 
that $Y_{L,S}=S$. Since $V$ is $H$-stable, contained in $S$ and $S$ is a single 
$H$-orbit, it follows that $V=S$ and $\ce$ is an $H$-equivariant local system on $S$ 
that is, $\ce\in\cs(S)$. Using 23.3(a), we see that $\ce$ is a cuspidal local system. 
The lemma is proved.

\subhead 30.3\endsubhead
For $J\sub\II$ such that $\e(J)=J$, let 
$$V_{J,D}=\{(P,gU_P);P\in\cp_J,gU_P\in N_DP/U_P\}.$$
Let $P_0\in\cp_J$. Since $\e(J)=J$, $N_DP_0$ is a connected component of $N_GP_0$ and 
$D_0=N_DP_0/U_{P_0}$ is a connected component of $N_GP_0/U_{P_0}$. Consider the diagram
$$D_0@<a<<G^0\T D_0@>b>>V_{J,D}$$
where $a(x,gU_{P_0})=gU_{P_0},b(x,gU_{P_0})=(xPx\i,xgx\i U_{xPx\i})$ with $x\in G^0$,
$g\in N_DP_0$. Then $a,b$ are smooth morphisms with connected fibres; more precisely, 
$b$ is a principal $P_0$-bundle where $P_0$ acts on $G^0\T D_0$ by
$p:(x,gU_{P_0})\m(xp\i,pgp\i U_{P_0})$. Let $A_0$ be a perverse sheaf on $D_0$
equivariant for the conjugation action of $P_0/U_{P_0}$. Then 
$a^\bst A_0=b^\bst A_0^\flat$ for a well defined perverse sheaf $A_0^\flat$ on 
$V_{J,D}$. 

\subhead 30.4\endsubhead
For any $J\sub J'\sub\II$ such that $\e(J)=J,\e(J')=J'$ let $V_{J,J',D}$ be the variety
consisting of all pairs $(P,gU_Q)$ where $P\in\cp_J,gU_Q\in N_DP/U_Q$ and $Q$ is the 
unique parabolic in $\cp_{J'}$ such that $P\sub Q$. We have a diagram
$$V_{J,D}@<c<<V_{J,J',D}@>d>>V_{J',D}$$
where $c(P,gU_Q)=(P,gU_P),d(P,gU_Q)=(Q,gU_Q)$. Define 
$$\tf_{J,J'}:\cd(V_{J,D})@>>>\cd(V_{J',D}),\te_{J,J'}:\cd(V_{J',D})@>>>\cd(V_{J,D})$$
by $\tf_{J,J'}A=d_!c^*A,\te_{J,J'}A'=c_!d^*A'$. Define
$$f_{J,J'}:\cd(V_{J,D})@>>>\cd(V_{J',D}),e_{J,J'}:\cd(V_{J',D})@>>>\cd(V_{J,D})$$
by $f_{J,J'}A=\tf_{J,J'}A[\a]=d_!c^\bst A,e_{J,J'}A'=\te_{J,J'}A'[\a](\a)$ where 
$\a=\dim\cp_J-\dim\cp_{J'}$. 

Let $P_0\in\cp_J,P'_0\in\cp_{J'}$ be such that $P_0\sub P'_0$. Let 
$D_0=N_DP_0/U_{P_0},D'_0=N_DP'_0/U_{P'_0}$. We have $\a=\dim U_{P_0}-\dim U_{P'_0}$. We
show:

(a) {\it If $A\in\cm(V_{J,D})$ is of the form $A=A_0^\flat$ where $A_0$ is a direct sum
of admissible simple perverse sheaves on $D_0$ then $A':=f_{J,J'}A$ is of the form 
$A'_0{}^\flat$ where $A'_0=\ind_{D_0}^{D'_0}A_0$ is a direct sum of admissible simple 
perverse sheaves on $D'_0$. In particular, $f_{J,J'}A\in\cm(V_{J',D})$.} 
\nl
We have a commutative diagram
$$\CD
D_0       @<r<<     V_1@>s>>        V_2@>t>>                           D'_0\\
@AaAA                   @AjAA             @AhAA                           @Aa'AA     \\
G^0\T D_0   @<1\T r<<G^0\T V_1@>1\T s>>G^0\T V_2            @>1\T t>>G^0\T D'_0\\
@VbVV                  @.            @VkVV        @Vb'VV\\
V_{J,D}@<c<<   {} @<<< V_{J,J',D}  @>d>> V_{J',D}
\endCD$$
Here 

$V_1=P'_0/U_{P'_0}\T N_DP_0/U_{P'_0}$, 

$V_2=\{(P,gU_{P'_0});P\in\cp_J,P\sub P'_0,g\in N_DP\}$,

$b$ is as in 30.3, $b'$ is the analogous map (with $P_0$ replaced by $P'_0$),

$a,j,h,a',r,t$ are given by the second projection,

$s(p'U_{P'_0},gU_{P'_0})=(p'P_0p'{}\i,p'gp'{}\i U_{P'_0})$ where 
$p'\in P'_0,g\in N_DP_0$,

$k(x,P,gU_{P'})=(xPx\i,xgx\i U_{xP'_0x\i})$ where $x\in G^0$, $(P,gU_{P'})\in V_2$.
\nl
All morphisms in this diagram (except $t,1\T t,d$) are smooth with connected fibres. 
Moreover, $s,b,k,b'$ are principal bundles with group $P_0/U_{P'_0},P_0,P'_0,P'_0$. We
may assume that $A_0,A$ are simple. There is a well defined simple perverse sheaf $A_1$
on $V_2$ such that $r^\bst A_0=s^\bst A_1$. We have $A'_0=t_!A_1$. Using the 
commutativity of the diagram above we see that $h^\bst A_1=k^\bst(c^\bst A)$. Since the
squares $(h,t,1\T t,a')$ and $(1\T t,b',k,d)$ are cartesian we have 
$(1\T t)_!h^\bst A_1=a'{}^\bst A_0$, 
$(1\T t)_!k^\bst(c^\bst A)=b'{}^\bst d_!(c^\bst A)=b'{}^\bst A'$. It follows that 
$a'{}^\bst A_0=b'{}^\bst A'$. From 27.2(d) we see that $A'_0=\ind_{D_0}^{D'_0}A_0$ is a
direct sum of admissible simple perverse sheaves on $D'_0$. Hence 
$a'{}^\bst A_0=b'{}^\bst A'$ is a direct sum of simple perverse sheaves on 
$G^0\T D'_0$. Hence $A'$ is a direct sum of simple perverse sheaves on $V_{J',D}$. This
proves (a).

We show:

(b) {\it Let $C_0$ be a $P'_0/U_{P'_0}$-equivariant simple perverse sheaf on $D'_0$ and
let $C=C_0^\flat$ be the corresponding simple perverse sheaf on $V_{J',D}$. Then for 
any $i\in\ZZ$ we have ${}^pH^i(e_{J,J'}C)=({}^pH^i(\res_{D'_0}^{D_0}C_0))^\flat$ 
(equality of perverse sheaves on $V_{J,D}$).}
\nl
We have a commutative diagram
$$\CD
D_0       @<f'<<     V_3@>f>>                           D'_0\\
@AaAA                   @AvAA                           @Aa'AA     \\
G^0\T D_0   @<1\T f'<<G^0\T V_3         @>1\T f>>G^0\T D'_0\\
@VbVV                 @VmVV        @Vb'VV\\
V_{J,D}@<c<<   V_{J,J',D}  @>d>> V_{J',D}
\endCD$$
where $a,b,a',b',c,d$ are as above, $V_3=N_DP_0/U_{P'_0}$, $f,f'$ are the obvious maps,
$v$ is the second projection and 
$$m(x,pU_{P'_0})=(xP_0x\i,xpx\i U_{xP'_0x\i}).$$
From this commutative diagram we see that $v^*f^*C_0[\dim U_{P'_0}]=m^*d^*C$ (we use 
that $b'{}^*C=a'{}^*C_0[\dim U_{P'_0}]$). Since the squares $(f',a,1\T f',v)$ and 
$(v,f,1\T f,a')$ are cartesian we have
$$a^*f'_!(f^*C_0)=(1\T f')_!v^*(f^*C_0), b^*c_! d^*C=(1\T f')_!m^*(d^*C)$$
hence $a^*f'_!(f^*C_0)[\dim U_{P'_0}]=b^*c_!d^*C$. Thus, 
$$a^*(\res_{D'_0}^{D_0}C_0)(-\a)[\dim U_{P'_0}]=b^*(e_{J,J'}C)[-\a](-\a).$$
Hence $a^*(\res_{D'_0}^{D_0}C_0)[\dim U_{P_0}]=b^*(e_{J,J'}C)$ and
$b^\bst(e_{J,J'}C)=a^\bst(\res_{D'_0}^{D_0}C_0)$,
$${}^pH^i(b^\bst(e_{J,J'}C))={}^pH^i(a^\bst(\res_{D'_0}^{D_0}C_0)).$$
Using this and \cite{\CS, I,(1.8.1)} we have
$$b^\bst({}^pH^i(e_{J,J'}C))=a^\bst({}^pH^i(\res_{D'_0}^{D_0}C_0))=
b^\bst(  ({}^pH^i(\res_{D'_0}^{D_0}C_0))^\flat).$$
Since $\b^\bst$ is fully faithful \cite{\CS, I,(1.8.3)}, we deduce the required
equality \lb ${}^pH^i(e_{J,J'}C)=({}^pH^i(\res_{D'_0}^{D_0}C_0))^\flat$.

\proclaim{Lemma 30.5} Let $A\in\cd(V_{J,D}),A'\in\cd(V_{J',D})$. We have
$$\Hom_{\cd(V_{J,D})}(e_{J,J'}A',A)=\Hom_{\cd(V_{J',D})}(A',f_{J,J'}A).$$
\endproclaim
Using the fact that $d$ is proper (hence $d_*=d_!$) and that $c$ is an affine space
bundle with fibres of dimension $\a$ (hence $c^!A=c^*A[2\a](\a)$), we have 
$$\align&\Hom(c_!d^*A'[\a](\a),A)=\Hom(d^*A'[\a](\a),c^!A)=
\Hom(d^*A'[\a](\a),c^*A[2\a](\a))\\&=\Hom(d^*A',c^*A[\a])=\Hom(A',d_*c^*A[\a])=
\Hom(A',d_!c^*A[\a]).\endalign$$
The lemma is proved.

\proclaim{Theorem 30.6} Let $A$ be a character sheaf on $D$. 

(a) Let $P_0$ be a parabolic of $G^0$ such that $N_DP_0\ne\em$ and let 
$D_0=N_DP_0/U_{P_0}$, a connected component of $N_GP_0/U_{P_0}$. Let $A_1$ be a 
character sheaf on $D_0$. Then $\ind_{D_0}^DA_1\in\cm(D)$.

(b) Let $P_0,D_0$ be as in (a). Then $\res_D^{D_0}A\in\cd(D_0)^{\le 0}$.

(c) Let $P_0,D_0$ be as in (a). Let $A_1$ be a character sheaf on $D_0$. Then 

$\Hom_{\cd(D_0)}(\res_D^{D_0}A,A_1)=\Hom_{\cd(D)}(A,\ind_{D_0}^DA_1)$.

(d) There exist $P_0,D_0$ as in (a) and a cuspidal character sheaf $A_1$ on $D_0$ such 
that $A$ is a direct summand of $\ind_{D_0}^DA_1$.

(e) $A$ is admissible.
\endproclaim
If $G^0=\{1\}$, the theorem is obvious. Assume now that $\dim G>0$ and that the theorem
is true when $G$ is replaced by a reductive group of dimension $<\dim G$. The proof of
the theorem for $G$ assuming this inductive assumption is given in 30.7-30.11.

\subhead 30.7\endsubhead
We show that 30.6(a) holds for $G$. If $P_0=G^0$ we have $D=D_0$, $\ind_{D_0}^DA_1=A_1$
and the result is obvious. Assume now that $P_0\ne G^0$. By 30.6(e) for 
$N_GP_0/U_{P_0}$, $A_1$ is admissible on $D_0$. Using 27.2(d) we see that 
$\ind_{D_0}^DA_1\in\cm(D)$, as required.

\subhead 30.8\endsubhead
We show that 30.6(b) holds for $G$. If $P_0=G^0$ we have $\res_D^{D_0}A=A\in\cm(D_0)$. 
Assume now that $P_0\ne G^0$. Let $J$ be such that $P_0\in\cp_J$. We identify 
$V_{\II,D}=D$ in the obvious way. We show:

(a) ${}^pH^i(e_{J,\II}A)=0$ for $i>0$.
\nl
Assume that this is not so; let $i$ be the largest integer such that
${}^pH^i(e_{J,\II}A)\ne 0$. Then $i>0$ and there exists a nonzero morphism 
$e_{J,\II}A@>>>{}^pH^i(e_{J,\II}A)[-i]$. Using Lemma 30.5 we deduce that
$$\Hom_{\cd(D)}(A,f_{J,\II}({}^pH^i(e_{J,\II}A)[-i]))\ne 0.$$
Using 30.4(b) we have
$$f_{J,\II}({}^pH^i(e_{J,\II}A))=f_{J,\II}(({}^pH^i(\res_D^{D_0}A))^\flat).$$
By 29.15, ${}^pH^i(\res_D^{D_0}A)$ is a finite direct sum of character sheaves on $D_0$
hence, by 30.6(e) for $N_GP_0/U_{P_0}$ it is a finite direct sum of admissible 
complexes on $D_0$. Using 30.4(a), we see that 
$C:=f_{J,\II}(({}^pH^i(\res_D^{D_0}A))^\flat)\in\cm(D)$. Thus we have 
$\Hom_{\cd(D)}(A,C[-i])\ne 0$ with $A,C\in\cm(D),i>0$. This contradicts
\cite{\CS, II,7.4}. Thus, (a) holds. 

Using (a) and 30.4(b) we see that for $i>0$ we have
$({}^pH^i(\res_D^{D_0}A))^\flat=0$ hence ${}^pH^i(\res_D^{D_0}A)=0$. It follows that
$\res_D^{D_0}A\in\cd(D_0)^{\le 0}$. Thus, 30.6(b) holds for $G$.

\subhead 30.9\endsubhead
We show that 30.6(c) holds for $G$. If $P_0=G^0$ the result is obvious. Assume now that
$P_0\ne G^0$. Let $J$ be as in 30.8. By 30.6(e) for $N_GP_0/U_{P_0}$, $A_1$ is 
admissible on $D_0$. Let $A_1^\flat$ be the simple perverse sheaf on $V_{J,D}$
corresponding to $A_1$ as in 30.3. From 30.4 we have 
$$\ind_{D_0}^DA_1=f_{J,\II}A_1^\flat,a^\bst e_{J,\II}A=b^\bst\res_D^{D_0}A,\tag a$$
where $a,b$ are as in 30.3. We have
$$\align&\Hom_{\cd(D_0)}(\res_D^{D_0}A,A_1)=
\Hom_{\cd(G^0\T D_0)}(b^\bst\res_D^{D_0}A,b^\bst A_1)\\&=
\Hom_{\cd(G^0\T D_0)}(a^\bst e_{J,\II}A,a^\bst A_1^\flat)=
\Hom_{\cd(V_{J,D})}(e_{J,\II}A,A_1^\flat)\\&=\Hom_{\cd(D)}(A,f_{J,\II}A_1^\flat)=
\Hom_{\cd(D)}(A,\ind_{D_0}^DA_1);\endalign$$
the first equality comes from \cite{\CS, I,(1.8.2)} and 30.6(b); the second equality 
comes from (a); the third equality comes from \cite{\CS, I,(1.8.2)} and 30.8(a); the 
fourth equality comes from 30.5; the fifth equality comes from (a). We see that 30.6(c)
holds for $G$.

\subhead 30.10\endsubhead
We show that 30.6(d) holds for $A$. If $A$ is cuspidal, we can take $P=G^0,A_1=A$ and 
the desired result holds. Thus, we may assume that $A$ is not cuspidal. Then there 
exist $P_0,D_0$ as in 30.6(a) such that $P\ne G^0$ and
$\res_D^{D_0}A[-1]\n\cd(D_0)^{\le 0}$. Then ${}^pH^i(\res_D^{D_0}A)\ne 0$ for some
$i\ge 0$. By 30.6(b) we have ${}^pH^j(\res_D^{D_0}A)=0$ for all $j>0$. It follows that 
${}^pH^0(\res_D^{D_0}A)\ne 0$ and there exists a non-zero morphism
$\res_D^{D_0}A@>>>{}^pH^0(\res_D^{D_0}A)$ in $\cd(D_0)$. Since ${}^pH^0(\res_D^{D_0}A)$
is a direct sum of character sheaves on $D_0$ (see 29.15) it follows that there exists 
a character sheaf $A_2$ on $D_0$ and a non-zero morphism $\res_D^{D_0}A@>>>A_2$ in 
$\cd(D_0)$. Using 30.6(c) it follows that there exists a non-zero morphism
$A@>>>\ind_{D_0}^DA_2$ in $\cd(D)$. By 30.6(a) this is a non-zero morphism in $\cm(D)$.
This must be injective since $A$ is simple. By the induction hypothesis, $A_2$ is a 
direct summand of a complex of the form $\ind_{D_1}^{D_0}A_3$ where $D_1=N_{D_0}Q/U_Q$,
$Q$ is a parabolic of $N_GP/U_P$ such that $N_{D_0}Q\ne\em$ and $A_3$ is a cuspidal 
character sheaf on $D_1$. By the induction hypothesis, $A_3$ is admissible. By the 
transitivity property 27.3(a) we have 
$\ind_{D_0}^D(\ind_{D_1}^{D_0}A_3)=\ind_{D_1}^DA_3$. Since $\ind_{D_0}^D$ commutes with
direct sums, we see that $\ind_{D_0}^DA_2$ is a direct summand of $\ind_{D_1}^DA_3$. 
Hence $A$ is isomorphic to a subobject of $\ind_{D_1}^DA_3$. From 27.2(d) we see that 
$\ind_{D_1}^DA_3$ is a semisimple perverse sheaf hence $A$ is a direct summand of it. 
Thus, 30.6(d) holds for $G$.

\subhead 30.11\endsubhead
We show that 30.6(e) holds for $A$. Assume first that $A$ is cuspidal. Then $A$ is 
admissible by Lemma 30.2, which is applicable, in view of 28.15(a) with $J=\II$. Next 
assume that $A$ is not cuspidal. Then, by 30.6(d) and its proof, we see that there 
exist $P_0,D_0$ as in 30.6(a) and a cuspidal character sheaf $A_1$ on $D_0$ such that 
$P_0\ne G^0$ and $A$ is a direct summand of $\ind_{D_0}^DA_1$. By the induction 
hypothesis, $A_1$ is admissible. Using 27.2(d) we see that $A$ is admissible. Thus 
30.6(e) holds. Theorem 30.6 is proved.

\proclaim{Corollary 30.12} Let $J\sub\II$ and let $X$ be a parabolic character sheaf on
$Z_{J,D}$ (see \S26). Then $X$ is admissible in the sense of 26.3.
\endproclaim
By \cite{\PCS, 4.13} we have $X=\hA$ where $\hA$ is obtained from some $\tt,C,A$ as in 
26.3 except that $A$ is a character sheaf on $C$ instead of being an admissible complex
on $C$. However, by 30.6(e), $A$ is automatically admissible on $C$ hence $\hA$ is 
admissible on $Z_{J,D}$, by the definition in 26.3.

\head 31. Character sheaves and Hecke algebras\endhead
\subhead 31.1\endsubhead
In this section we show that the restriction functor studied in \S29 takes a character
sheaf to a direct sum of character sheaves (Theorem 31.14). In the connected case this
result appeared in \cite{\CS, I,\S6} with a proof based on a connection of character 
sheaves with Hecke algebras. The present proof in the general case is an extension of
the proof in \cite{\CS, I,\S6}, taking also into account the approach given later by 
Mars and Springer \cite{\MS, \S9}. 

\subhead 31.2\endsubhead
Until the end of 31.13 we fix $n\in\NN_{\kk}^*$. We write $\fs_n$ instead of 
$\fs_n(\TT)$. Let $\ufs_n$ be the set of isomorphism classes of objects in $\fs_n$. We
have canonically 
$$\ufs_n=\Hom(\mu_n(\TT),\bbq^*)$$ 
(notation of 5.3). Thus $\sh(\ufs_n)=n^{\dim\TT}<\iy$. Now $\WW^\bul$ acts on $\ufs_n$ 
by $a:\l\m a\l$ where $\l$ is the isomorphism class of $\cl\in\fs_n$ and $a\l$ is the 
isomorphism class of $(a\i)^*\cl\in\fs_n$. For $\l\in\ufs_n$ we set $\WW_\l=\WW_\cl$ 
(see 28.3) where $\l$ is the isomorphism class of $\cl\in\fs_n$ (this is independent of
the choice of $\cl$).

Let $\ca=\ZZ[v,v\i]$ where $v$ is an indeterminate. We now introduce an associative 
$\ca$-algebra $H_n$ closely related to the algebra of double cosets of a finite 
Chevalley group with respect to a maximal unipotent subgroup, studied in \cite{\YO}. We
define $H_n$ by the generators $T_w(w\in\WW),1_\l(\l\in\ufs_n)$ and the relations

$1_\l1_{\l'}=\d_{\l,\l'}1_\l$ for $\l,\l'\in\ufs_n$,

$T_wT_{w'}=T_{ww'}$ for $w,w'\in\WW$ with $l(ww')=l(w)+l(w')$,

$T_w1_\l=1_{w\l}T_w$ for $w\in\WW,\l\in\ufs_n$,

$T_s^2=v^2T_1+(v^2-1)\sum_{\l;s\in\WW_\l}T_s1_\l$ for $s\in\II$,

$T_1=\sum_\l1_\l$.
\nl
Note that $T_1=\sum_\l1_\l$ is the unit element of $H_n$ and that

(a) {\it $\{T_w1_\l;w\in\WW,\l\in\ufs_n\}$ is an $\ca$-basis of $H_n$.}
\nl
In the case where $n=1$, $H_n$ is just the Hecke algebra attached to $\WW$ and the 
proof of (a) is standard (see for example \cite{\HA, 3.3}). The proof in the general 
case is quite similar. Consider the free $\ca$-module $M$ with basis 
$$\{[\l',w,\l];\l,\l'\in\ufs_n,w\in\WW,w\l=\l'\}.$$ 
If (a) is true, then $M$ may be identified with $H_n$ so that $[\l',w,\l]\in M$ 
corresponds to $T_w1_\l=1_{\l'}T_w\in H_n$; hence $M$ is naturally an $(H_n,H_n)$ 
bimodule. Conversely, if we can make $M$ naturally into an $(H_n,H_n)$ bimodule then 
(a) can be easily deduced. For any $s\in\II$ we define $\ca$-linear maps $m\m T_sm$ and
$m\m mT_s$ of $M$ into itself by

$T_s[\l',w,\l]=[s\l',sw,\l]$ if $l(sw)=l(w)+1$,

$T_s[\l',w,\l]=v^2[s\l',sw,\l]+(v^2-1)[s\l',w,\l]$ if $l(sw)=l(w)-1,s\in\WW_{\l'}$,

$T_s[\l',w,\l]=v^2[s\l',sw,\l]$ if $l(sw)=l(w)-1,s\n\WW_{\l'}$,

$[\l',w,\l]T_s=[\l',ws,s\l]$ if $l(ws)=l(w)+1$, 

$[\l',w,\l]T_s=v^2[\l',ws,s\l]+(v^2-1)[\l',w,s\l]$ if $l(ws)=l(w)-1,s\in\WW_\l$,

$[\l',w,\l]T_s=v^2[\l',ws,s\l]$ if $l(ws)=l(w)-1,s\n\WW_\l$.
\nl
For any $\til\in\ufs_n$ we define $\ca$-linear maps $m\m 1_{\til}m,m\m m1_{\til}$ of 
$M$ into itself by

$1_{\til}[\l',w,\l]=\d_{\til,\l'}[\l',w,\l],[\l',w,\l]1_{\til}=\d_{\til,\l}[\l',w,\l]$.
\nl 
One shows that this defines an $(H_n,H_n)$ bimodule structure on $M$. We omit further
details.

{\it Remark.} In \cite{\MS, 3.3.1} an algebra structure on $M$ (with $\l$ restricted in
a fixed $W$-orbit) is considered which is similar to the one above, coming from $H_n$, 
but differs from it in the following way: for $s\in\II$ and $\l$ such that $s\n\WW_\l$,
$[s\l,s,\l]^2$ is equal, in our definition, to $v^2[s\l,1,\l]$, while in the definition
of \cite{\MS} it equals $[s\l,1,\l]$.

\subhead 31.3\endsubhead
We return to the general case. For $s\in\II$, $T_s$ is invertible in $H_n$; we have
$$T_s\i=v^{-2}T_s+(v^{-2}-1)\sum_{\l;s\in\WW_\l}1_\l.$$
Moreover,
$$T_s\i T_s\i=v^{-2}+(v^{-2}-1)\sum_{\l;s\in\WW_\l}T_s\i 1_\l.$$
It follows that that $T_w$ is invertible in $H_n$ for any $w\in\WW$ and that 
$$T_{y\i}\i T_{w\i}\i=T_{(yw)\i}\i\text{ if }y,w\in\WW,l(yw)=l(y)+l(w).$$
For any $w\in\WW,\l\in\ufs_n$ we have $T_{w\i}\i 1_\l=1_{w\l}T_{w\i}\i$. Hence there is
a unique ring homomorphism $\bar{}:H_n@>>>H_n$ such that $\ov{T_w}=T_{w\i}\i$ for all 
$w\in\WW$, $\ov{v^m1_\l}=v^{-m}1_\l$ for all $\l$ and all $m\in\ZZ$. Note that

(a) {\it $\bar{}:H_n@>>>H_n$ has square $1$. In particular, it is an isomorphism of
rings.}
\nl
Indeed, the generators $v,v\i,T_w(w\in\WW)$ and $1_\l(\l\in\ufs_n)$ are mapped to 
themselves by the square of $\bar{}$.

\subhead 31.4\endsubhead
In the remainder of this section we fix a connected component $D$ of $G$. We write 
$\e:\WW@>>>\WW$ instead of $\e_D:\WW@>>>\WW$. Let $\o$ be the order of 
$\io_D:\TT@>>>\TT$. The assignment $1_\l\m 1_{\uD\l},T_w\m T_{\e(w)}$ defines an 
automorphism of the algebra $H_n$ denoted by $h\m\fa_D(h)$. We have $\fa_D^\o=1$. Let 
$H_n\la[D]\ra$ be the free left $H_n$-module with basis $\{[D]^i;i\in\ZZ/\o\ZZ\}$. We 
regard $H_n\la[D]\ra$ as an associative $\ca$-algebra with unit $1=[D]^0$ such that 
$$(h[D]^i)(h'[D]^{i'})=h\fa_D^i(h')[D]^{i+i'}$$
for $h,h'\in H_n,i,i'\in\ZZ/\o\ZZ$. We write $[D]$ instead of $D^1$. We have
$\fa_D(h)=[D]h[D]\i$ in $H_n\la[D]\ra$. Define a group involution 
$\bar{}:H_n\la[D]\ra@>>>H_n\la[D]\ra$ by $\ov{h[D]^i}=\ov{h}[D]^i$ where $h\in H_n$,
$i\in\ZZ/\o\ZZ$ and $\bar{}:H_n@>>>H_n$ is as in 31.3. Then
$\bar{}:H_n\la[D]\ra@>>>H_n\la[D]\ra$ is a ring involution (we use the fact that 
$\bar{}:H_n@>>>H_n$ commutes with $\fa_D$).

\subhead 31.5\endsubhead
For $s\in\II\cup\{1\},\l\in\ufs_n$ we set 
$$C^s_\l=(T_s+u)1_\l=1_{s\l}(T_s+u)\in H_n,$$
where $u=1$ if $s\in\WW_\l\cap\II$ and $u=0$ if $s\n\WW_\l$. We have
$$\ov{C^s_\l}=v^{-2}C^s_\l\text{ if }s\in\II,\qua \ov{C^1_\l}=C^1_\l.\tag a$$ 
If $\ss=(s_1,s_2,\do,s_r)$ is a sequence in $\II\cup\{1\}$ and $\l\in\ufs_n$, we set
$$C^\ss_\l=(C^{s_1}_{s_2\do s_r\l})(C^{s_2}_{s_3\do s_r\l})\do(C^{s_r}_\l)\in H_n.$$
It is easy to check that 

(b) {\it the $\ca$-module of $H_n$ is generated by the elements $C^\ss_\l$ with 
$\ss,\l$ as above.}

\subhead 31.6\endsubhead
If $A$ is a simple perverse sheaf on an algebraic variety $V$ and $K$ is a perverse 
sheaf on $V$ we write $(A:K)$ for the multiplicity of $A$ in a Jordan-H\"older series 
of $K$. 

Until the end of 31.9, $A$ denotes a fixed character sheaf on $D$. For

(a) $\l\in\ufs_n$ and a sequence $\ss=(s_1,s_2,\do,s_r)$ in $\II\cup\{1\}$ with 
$s_1s_2\do s_r\uD\l=\l$
\nl
we set 
$$\g_\l(\ss)=\sum_{j\in\ZZ}(-v)^j(A:{}^pH^j(\bK^{\ss,\cl}_D))\in\ca$$ 
where $\cl\in\fs_n$ is in the isomorphism class $\l$.

\proclaim{Proposition 31.7} There is a unique $\ca$-linear map $\z^A:H_n[D]@>>>\ca$ 
such that for any $\l\in\ufs_n$ and any sequence $\ss=(s_1,s_2,\do,s_r)$ in 
$\II\cup\{1\}$ we have 
$$\z^A(C^\ss_{\uD\l}[D])=v^{-\dim G}\g_\l(\ss)\text{ if }s_1s_2\do s_r\uD\l=\l,$$
$$\z^A(C^\ss_{\uD\l}[D])=0\text{ if }s_1s_2\do s_r\uD\l\ne\l.$$
\endproclaim
Let $X$ be the free $\ca$-module with basis $[\l',\ss,\l]$ where 
$\ss=(s_1,s_2,\do,s_r)$ is a sequence in $\II\cup\{1\}$ and $\l,\l'\in\ufs_n$ are such
that $s_1s_2\do s_r\uD\l=\l'$. Define a $\ca$-linear map $c:X@>>>H_n[D]$ by 
$c[\l',\ss,\l]=C^\ss_{\uD\l}[D]$. Define a $\ca$-linear map $c':X@>>>\ca$ by 
$c'[\l',\ss,\l]=v^{-\dim G}\g_\l(\ss)$ if $\l'=\l$, $c'[\l',\ss,\l]=0$ if $\l'\ne\l$. 
Since $c$ is surjective (see 31.5(b)) we see that it suffices to prove the following
statement.

(a) {\it If $\x\in\ker(c)$ then $\x\in\ker(c')$.}
\nl
By a standard argument \cite{\BBD, \S6}, to prove (a), we may assume that

(b) {\it $\kk$ is an algebraic closure of a finite field $\FF_q$, that $G$ has a fixed 
$\FF_q$-structure with Frobenius map $F:G@>>>G$ which induces the identity map on $\WW$
and on $G/G^0$ and the map $t\m t^q$ on $\TT$ and that $F^*\cl\cong\cl$ for any
$\cl\in\fs_n$.}
\nl
(We use that $\ufs_n$ is a finite set.) Then each $\cl$ as in (a) may be regarded as 
pure of weight $0$ and each $K^{\cl,w}_D,K^{\cl,\ww}_D,\bK^{\cl,\ss}_D$ (as in 28.13) 
may be regarded as a mixed complex on $D$. For any mixed complex $K$ on $D$ we denote 
by ${}^pH^i_j(K)$ the $j$-th subquotient of the weight filtration of ${}^pH^i(K)$ so 
that ${}^pH^i_j(K)$ is a pure perverse sheaf of weight $j$; we set
$$\c^A_v(K)=\sum_{i,j}(-1)^i(A:{}^pH^i_j(K))v^j\in\ca.$$
Define a $\ca$-linear map $\tiz^A:H_n[D]@>>>\ca$ by
$$\tiz^A(T_w1_{\uD\l}[D])=v^{-\dim G}\c^A_v(K^{w,\cl}_D),\text{ if }w\uD\l=\l,$$
$$\tiz^A(T_w1_{\uD\l}[D])=0\text{ if }w\uD\l\ne\l$$
where $\cl\in\fs_n$ is in the isomorphism class $\l$. 

Let $\l,\cl,\ss,r$ be as in 31.6(a); we will show that
$$\c^A_v(K^{\ss,\cl}_D)=v^{\dim G}\tiz^A(T_{s_1}T_{s_2}\do T_{s_r}1_{\uD\l}[D]),\tag c
$$
$$\c^A_v(\bK^{\ss,\cl}_D)=v^{\dim G}\tiz^A(C^\ss_{\uD\l})[D]).\tag d$$
We prove (c) by induction on $r$. We may assume that all $s_i$ are in $\II$. If $r=0$ 
we have $K^{\ss,\cl}_D=K^{1,\cl}_D$ and the result is clear; we have
$$\c^A_v(K^{1,\cl}_D)=v^{\dim G}\tiz^A(1_{\uD\l}[D]).$$
Assume that $r\ge 1$. Assume first that $l(s_1s_2\do s_r)=r$. Using 28.13(a) repeatedly
we have $K^{\ss,\cl}_D=K^{w,\cl}_D$ where $w=s_1s_2\do s_r$. Then the right hand side 
of (c) equals $v^{\dim G}\tiz^A(T_w1_{\uD\l})$; the result follows. Assume next that
$l(s_1s_2\do s_r)<r$. We can find $j\in[2,r]$ such that $s_j\do s_{r-1}s_r$ is a 
reduced expression in $\WW$ and $s_{j-1}s_j\do s_{r-1}s_r$ is not a reduced expression.
We can find $s'_j,\do,s'_{r-1},s'_r$ in $\II$ such that 
$s'_j\do s'_{r-1}s'_r=s_j\do s_{r-1}s_r$ and $s'_j=s_{j-1}$. Let 
$$\uu=(s_1,s_2,\do,s_{j-1},s'_j,\do,s'_{r-1},s'_r).$$
 As in the proof of the implication
(v)$\imp$(i) in 28.13 we see that $K^{\ss,\cl}_D=K^{\uu,\cl}_D$. Hence
$$\c^A_v(K^{\ss,\cl}_D)=\c^A_v(K^{\uu,\cl}_D).$$
Moreover, both the right hand side of (c) and the analogous expression obtained by 
replacing $\ss$ by $\uu$ are equal to 
$$v^{\dim G}\tiz^A((T_{s_1}\do T_{s_{j-1}}T_{s_js_{j+1}\do s_r}1_{\uD\l})[D]).$$
Hence to prove the lemma for $\ss$ it suffices to prove it for $\uu$. Thus we are 
reduced to the case where $s_j=s_{j-1}$. In this case we will use the notation in
28.13(d),(e),(f) with $J=\II$.

If $j\n\cj_\ss$, then from 28.13(f) we have 
$\c^A_v(K^{\ss,\cl}_D)=v^2\c^A_v(K^{\ss'',\cl}_D)$. In this case we have
$$T_{s_{j-1}}T_{s_j}1_{s_{j+1}\do s_r\uD\l}=v^21_{s_{j+1}\do s_r\uD\l}$$
since $s_j\n\WW_{s_{j+1}\do s_r\uD\l}$. Hence the right hand side of (c) is equal to
$$v^2v^{\dim G}\tiz^A(T_{s_1}\do T_{s_{j-2}}T_{s_{j+1}}T_{s_{j+2}}T_{s_r}1_{\uD\l}[D])
$$
which by the induction hypothesis is equal to $v^2\c^A_v(K^{\ss'',\cl}_D)$. Thus, (c) 
holds in this case.

If $j\in\cj_\ss$, then from 28.13(d),(e) we have 
$$\c^A_v(K^{\ss,\cl}_D)=\c^A_v(\p_{1!}\tcl)+v^2\c^A_v(K^{\ss'',\cl}_D),
\qua v^2\c^A_v(K^{\ss',\cl}_D)=\c^A_v(\p_{1!}\tcl)+\c^A_v(K^{\ss',\cl}_D).$$
(Indeed since the weight filtrations are strictly compatible with morphisms
\cite{\BBD, 5.3.5}, the exact sequences 28.13(d),(e) remain exact when each ${}^pH^i$
is replaced by its pure subquotient of a fixed weight.) It follows that
$$\c^A_v(K^{\ss,\cl}_D)=v^2\c^A_v(K^{\ss'',\cl}_D)+(v^2-1)\c^A_v(K^{\ss',\cl}_D).$$
Using the induction hypothesis for $\ss'',\ss'$ we see that
$$\align&v^2\c^A_v(K^{\ss'',\cl}_D)+(v^2-1)\c^A_v(K^{\ss',\cl}_D)=v^2v^{\dim G}
\tiz^A(T_{s_1}\do T_{s_{j-2}}T_{s_{j+1}}T_{s_r}1_{\uD\l}[D])\\&+(v^2-1)v^{\dim G}
\tiz^A(T_{s_1}\do T_{s_{j-1}}T_{s_{j+1}}\do T_{s_r}1_{\uD\l}[D]).\endalign$$
Substituting here
$$v^2T_11_{s_{j+1}\do s_r\uD\l})+(v^2-1)T_{s_{j-1}}1_{s_{j+1}\do s_r\uD\l})
=T_{s_{j-1}}T_{s_j}1_{s_{j+1}\do s_r\uD\l})$$
which holds since $s_{j-1}=s_j\in\WW_{s_{j+1}\do s_r\uD\l}$, we see that 
$v^2\c^A_v(K^{\ss'',\cl}_D)+(v^2-1)\c^A_v(K^{\ss',\cl}_D)$ is equal to the right hand 
side of (c). Thus (c) holds.

We prove (d). We will use the notation in 28.13(b) with $J=\II$. Using 28.13(b) we get
$$\c^A_v(\bpi_\ss{}_!f^i_!(f^i)^*\bcl)=\c^A_v(\bpi_\ss{}_!f^{i+1}_!(f^{i+1})^*\bcl)
+\sum_{\cj\sub\cj_\ss;|\cj|=i}\c^A_v(K^{\ss_\cj,\cl}_D)$$
for any $i$. Summing these equalities over all $i\ge 0$ we find
$$\c^A_v(\bK^{\ss,\cl}_D)=\sum_{\cj\sub\cj_\ss}\c^A_v(K^{\ss_\cj,\cl}_D).$$
We now use (c) for each $\ss_\cj$ in the last sum. We see that 
$\c^A_v(\bK^{\ss,\cl}_D)$ is a sum of $2^k$ terms (each term is $v^{\dim G}$ times a 
product of basis elements of $H_n$ times $[D]$) where $k=|\cj_\ss|$. Clearly, the right
hand side of (d) is the sum of the same $2^k$ terms. This proves (d). 

As in the proof of 29.14, $\bK^{\ss,\cl}_D$ is pure of weight $0$ hence
${}^pH^i_j(\bK^{\ss,\cl}_D)$ equals ${}^pH^i(\bK^{\ss,\cl}_D)$ if $i=j$ and equals $0$
if $i\ne j$. It follows that
$$\c^A_v(\bK^{\ss,\cl}_D)=\g_\l(\ss)$$
hence (d) implies $c'[\l,\ss,\l]=\tiz^A(C^\ss_{\uD\l})[D])=\tiz^A(c[\l,\ss,\l])$.
Clearly, $c'[\l',\ss,\l]=\tiz^A(c[\l',\ss,\l])$ if $\l\ne\l'$ (both sides are $0$). We
see that $c'=\tiz^A\circ c$. Hence (a) holds. The proposition is proved.

\proclaim{Lemma 31.8}For any $h,h'\in H_n$ we have $\z^A(hh'[D])=\z^A(h'[D]h)$.
\endproclaim
By 31.5(b), we may assume that $h=a_1\do a_{p-1}$, $h'=a_pa_{p+1}\do a_r$ where
$$a_i=C^{s_i}_{s_{i+1}\do s_{p-1}\uD\til}\text{ for }i\in[1,p-1],$$
$$a_i=C^{s_i}_{s_{i+1}\do s_r\uD\l}\text{ for }i\in[p,r];$$
here $\ss=(s_1,s_2,\do,s_r)$ is a sequence in $\II$, $1\le p\le r$ and 
$\l,\til\in\ufs_n$. If $s_p\do s_r\uD\l\ne\uD\til$ or $s_1\do s_{p-1}\uD\til\ne\l$ then
$\z^A(hh'[D])=0$ and $\z^A(h'[D]h)=0$. Thus we may assume that 
$\uD\til=s_p\do s_r\uD\l$ and $s_1\do s_{p-1}\uD\til=\l$. Then we have 
$$a_i=C^{s_i}_{s_{i+1}\do s_r\uD\l}\text{ for }i\in[1,r]$$ 
and $s_1s_2\do s_r\uD\l=\l$. Hence, setting $\ss=(s_1,s_2,\do,s_r)$, we see that 
$\bK^{\ss,\cl}_D$ is defined (with $\cl\in\fs_n$ in the isomorphism class $\l$). Let 
$$\ss'=(s_p,s_{p+1},\do,s_r,\e(s_1),\do,\e(s_{p-1})),$$
let $\l'=s_{p-1}\do s_1\l$ and let $\cl'\in\fs_n$ be in the isomorphism class $\l'$. 
Then 
$$s_ps_{p+1}\do s_r\e(s_1)\do\e(s_{p-1})\uD\l'=\l'$$ 
hence $\bK^{\ss',\cl'}_D$ is defined. Using 28.16(a) $p-1$ times, we have 
$\bK^{\ss,\cl}_D=\bK^{\ss',\cl'}_D$. Hence $\g_\l(\ss)=\g_{\l'}(\ss')$. By 31.7, we 
have $\z^A(a_1a_2\do a_r[D])=\z^A(a'_1a'_2\do a'_r[D])$ with $a_i$ as above and
$$a'_i=C^{s_{p+i-1}}_{s_{p+i}\do s_r\e(s_1\do s_{p-1})\uD\l'}
=C^{s_{p+i-1}}_{s_{p+i}\do s_r\uD\l}=a_{p+i-1}$$
for $i\in[1,r-p+1]$,   
$$a'_i=C^{\e(s_{i-r+p-1})}_{\e(s_{i-r+p}\do s_{p-1})\uD\l'}
=C^{\e(s_{i-r+p-1})}_{\uD s_{i-r+p}\do s_r\uD \l}=[D]a_{i-r+p-1}[D]\i$$
for $i\in[r-p+2,r]$. We have therefore $\z^A(hh'[D])=\z^A(h'[D]h)$. The lemma is  
proved.

\proclaim{Lemma 31.9}For any $h\in H_n$ we have $\z^A(\ov{h}[D])=\ov{\z^A(h[D])}$. Here
$\bar{}:\ca@>>>\ca$ is the ring homomorphism such that $\ov{v^m}=v^{-m}$ for all $m$.
\endproclaim
Using 31.5(b) we may assume that $h=C^\ss_\l$ where $\ss=(s_1,s_2,\do,s_r)$ is a 
sequence in $\II$ and $\l\in\ufs_n$. Using 31.5(a) we see that $\ov{h}=v^{-2r}h$. Hence
it suffices to show that $\ov{\z^A(h[D])}=v^{-2r}\z^A(h[D])$. If 
$s_1s_2\do s_r\uD\l\ne\l$ then $\z^A(h[D])=0$ and the result is obvious. Thus we may
assume that $s_1s_2\do s_r\uD\l=\l$ so that $\bK^{\ss,\cl}_D$ is defined (with 
$\cl\in\fs_n$ in the isomorphism class $\l$). By Lemma 31.7 we have 
$\z^A(h[D])=v^{-\dim G}\g_\l(\ss)$. Hence it suffices to show that
$$\ov{\g_\l(\ss)}=v^{-2m}\g_\l(\ss)$$
where $m=r+\dim G$. Using 28.17(a) we have
$$\align&\ov{\g_\l(\ss)}=\sum_j(-v)^{-j}(A:{}^pH^j(\bK^{\ss,\cl}_D))
=\sum_j(-v)^{-j}(A:{}^pH^{2m-j}(\bK^{\ss,\cl}_D))\\&
=\sum_j(-v)^{j-2m}(A:{}^pH^j(\bK^{\ss,\cl}_D))=v^{-2m}\g_\l(\ss).\endalign$$
The lemma is proved.

\proclaim{Lemma 31.10}Let $J\sub\II,y\in{}^J\WW,s\in\II,\l'\in\ufs_n$. We have
$$T_yC^s_{\l'}=v^{2\d(y)}\sum_{y_1}(C^{t_1}_{y_1\l'})T_{y_1};\tag a$$
here $\d(y)$ is $1$ if $ys<y,ys\in{}^J\WW$ and is $0$ otherwise; the sum is taken over 
all $y_1\in{}^J\WW\cap\{y,ys\}$ such that
$$ys\n\WW_Jy,y_1=y\imp s\in\WW_{\l'};$$
$t_1\in J\cup\{1\}$ is defined by $ys=t_1y$ if $ys\in\WW_Jy$ and $t_1=1$ if 
$ys\n\WW_Jy$.
\endproclaim  
If $s\n\WW_{\l'},ys>y$ then both sides of (a) are equal to $T_{ys}1_{\l'}$. If 
$s\n\WW_{\l'},ys<y$ then both sides of (a) are equal to $v^2T_{ys}1_{\l'}$. If 
$s\in\WW_{\l'},ys>y$ then both sides of (a) are equal to $(T_{ys}+T_y)1_{\l'}$. If 
$s\in\WW_{\l'},ys<y$ then both sides of (a) are equal to $v^2(T_{ys}+T_y)1_{\l'}$. The 
lemma is proved.

\proclaim{Lemma 31.11}Let $J\sub\II,y\in{}^J\WW,\l'\in\ufs_n$ and let
$\ss=(s_1,s_2,\do,s_r)$ be a sequence in $\II$. We have
$$T_yC^\ss_{\l'}=\sum_\yy v^{2\d(\yy)}(C^\tt_{y_r\l'})T_{y_r};$$
here the sum is taken over all sequences $\yy=(y_0,y_1,\do,y_r)$ in ${}^J\WW$ such that
$y=y_0$ and $y_i\in\{y_{i-1},y_{i-1}s_i\}$ for $i\in[1,r]$, $\tt=(t_1,t_2,\do,t_r)$ is 
the sequence in $J\cup\{1\}$ defined by $y_{i-1}s_i=t_iy_{i-1}$ if 
$y_{i-1}s_i\in\WW_Jy_{i-1}$ and $t_i=1$ if $y_{i-1}s_i\n\WW_Jy_{i-1}$; these are
subject to the requirement 
$$i\in[1,r],t_i=1,y_{i-1}=y_i\imp s_i\in\WW_{s_{i+1}\do s_r\l'};$$
moreover,
$$\d(\yy)=\sh(i\in[1,r];y_{i-1}s_i<y_{i-1},y_{i-1}s_i\in{}^J\WW).\tag a$$
\endproclaim
This follows by applying $r$ times Lemma 31.10.

\subhead 31.12\endsubhead
Until the end of 31.14 we fix $D,P,L,G',D'$ as in 29.1. Let $J\sub\II$ be such that 
$P\in\cp_J$. Let $H_{J,n}$ be the $\ca$-algebra defined in terms of $L$ in the same way
as $H_n$ was defined in 31.2 in terms of $G^0$. Since the Weyl group of $L$ is 
naturally the subgroup $\WW_J$ of $\WW$ and the canonical torus of $L$ may be 
identified with $\TT$ as in 29.1, we may identify $H_{J,n}$ with the subalgebra of 
$H_n$ generated as an $\ca$-submodule by 

(a) $\{T_w1_\l;w\in\WW_J,\l\in\ufs_n\}$.
\nl
Then $H_n$ is naturally a left $H_{J,n}$-module (using left multiplication). This 
$H_{J,n}$-module is free with basis $\{T_y;y\in{}^J\WW\}$. (The elements $bT_y$ where 
$b$ runs through the set (a) and $y\in{}^J\WW$, form the $\ca$-basis 31.2(a) of $H_n$.)
Applying to this basis the ring involution $\bar{}:H_n@>>>H_n$ which restricts to a 
ring involution of $H_{J,n}$ we deduce that $\{T_{y\i}\i;y\in{}^J\WW\}$ is a basis of
the left $H_{J,n}$-module $H_n$.

\proclaim{Lemma 31.13} Let $\l\in\ufs_n$ and let $\ss=(s_1,s_2,\do,s_r)$ be a sequence
in $\II$ such that $s_1s_2\do s_r\uD\l=\l$. Let $\cl\in\fs_n$ be in the isomorphism
class $\l$. Let $m=r+\dim G$. For any $j\in\ZZ$ we have
$${}^pH^j(\res_D^{D'}(\bK^{\ss,\cl}_D))\cong{}^pH^{2m-j}(\res_D^{D'}(\bK^{\ss,\cl}_D)).
$$
\endproclaim
Let $\Xi$ be the set of all pairs $(\yy,\tt)$ where $\yy=(y_0,y_1,\do,y_r)$ is a 
sequence in ${}^J\WW$ such that $y_r=\e(y_0)$ and $y_i\in\{y_{i-1},y_{i-1}s_i\}$ for
$i\in[1,r]$, $\tt=(t_1,t_2,\do,t_r)$ is the sequence in $J\cup\{1\}$ defined by 
$y_{i-1}s_i=t_iy_{i-1}$ if $y_{i-1}s_i\in\WW_Jy_{i-1}$ and $t_i=1$ if 
$y_{i-1}s_i\n\WW_Jy_{i-1}$; these are subject to the requirement 
$$i\in[1,r],t_i=1,y_{i-1}=y_i\imp s_i\in\WW_{s_{i+1}\do s_r\uD\l}.$$
Using 29.14 we see that it suffices to show that for any $j$ we have
$$\op_{(\yy,\tt)\in\Xi}{}^pH^{j-2\dd(\yy)}(\bK^{\tt,{}^{y_0}\cl}_{D'}))\cong
\op_{(\yy,\tt)\in\Xi}{}^pH^{2m-j-2\dd(\yy)}(\bK^{\tt,{}^{y_0}\cl}_{D'})).$$
Here both sides are semisimple complexes (see 28.12); hence it suffices to show that 
for any character sheaf $A'$ on $D'$ we have
$$\align&
\sum_{(\yy,\tt)\in\Xi;j}(-v)^j(A':{}^pH^{j-2\dd(\yy)}(\bK^{\tt,{}^{y_0}\cl}_{D'}))\\&=
\sum_{(\yy,\tt)\in\Xi;j}(-v)^j(A':{}^pH^{2m-j-2\dd(\yy)}(\bK^{\tt,{}^{y_0}\cl}_{D'})).
\endalign$$
or equivalently
$$\align& v^{2\dd(\yy)}\sum_{(\yy,\tt)\in\Xi;j}(-v)^j(A':{}^pH^j
(\bK^{\tt,{}^{y_0}\cl}_{D'}))\\&=v^{2m-2\dd(\yy)}
\sum_{(\yy,\tt)\in\Xi;j}(-v)^{-j}(A':{}^pH^j(\bK^{\tt,{}^{y_0}\cl}_{D'})).\endalign$$
Using 31.7 for $G',A',\tt$ instead of $G,A,\ss$, we see that it suffices to show 
$$\sum_{(\yy,\tt)\in\Xi}\z^{A'}(C^\tt_{\uD y_0\l}[D])v^{\dim L}v^{2\dd(\yy)}
=\sum_{(\yy,\tt)\in\Xi}\ov{\z^{A'}(C^\tt_{\uD y_0\l}[D])v^{\dim L}v^{2\dd(\yy)}}v^{2m}.
$$
Here $\z^{A'}:H_{J,n}[D]@>>>\ca$ is defined as in 3.7 for $G',A'$ instead of $G,A$. 
Here we substitute $\dd(\yy)=\d(y)+\dim U_P$ with $\d(\yy)$ as in 31.11(a), and use 
$\uD y_0\l=y_r\uD\l$ and $\dim L+2\dim U_P=\dim G$; we see that it suffices to show
$$\z^{A'}(\Psi[D])=\ov{\z^{A'}(\Psi[D])}v^{2r}$$
where $\Psi=\sum_{(\yy,\tt)\in\Xi}C^\tt_{y_r\uD\l}v^{2\d(\yy)}\in H_{J,n}$. We write 
the matrix of right multiplication by $C^\ss_{\uD\l}$ in $H_n$ (an $H_{J,n}$-linear 
map) in the $H_{J,n}$-basis $\{T_y;y\in{}^J\WW\}$:
$$T_yC^\ss_{\uD\l}=\sum_{y'\in{}^J\WW}a_{y,y'}T_{y'}\tag a$$
where $y\in{}^J\WW$ and $a_{y,y'}\in H_{J,n}$. Using Lemma 31.11 (with $\l'=\uD\l$) we 
see that $\sum_ya_{y,\e(y)}=\Psi$ where $y$ runs over ${}^J\WW$. Hence it suffices to 
show that
$$\z^{A'}(\sum_ya_{y,\e(y)}[D])=\ov{\z^{A'}(\sum_ya_{y,\e(y)}[D])}v^{2r}.$$
Using Lemma 31.9 (for $G',A'$ instead of $G,A$) we see that
$$\ov{\z^{A'}(\sum_ya_{y,\e(y)}[D])}=\z^{A'}(\sum_y\ov{a_{y,\e(y)}}[D]).$$
Hence it suffices to show that
$$\sum_y\z^{A'}(a_{y,\e(y)}[D])=\sum_y\z^{A'}(\ov{a_{y,\e(y)}}[D])v^{2r}.$$
Applying $\bar{}:H_n@>>>H_n$ to the equality (a) and using 
$\ov{C^\ss_{\uD\l}}=v^{-2r}C^\ss_{\uD\l}$ (see the proof of Lemma 31.9) we see that
$$v^{-2r}T_{y\i}\i C^\ss_{\uD\l}=\sum_{y'\in{}^J\WW}\ov{a_{y,y'}}T_{y'{}\i}\i.$$
Since $\{T_y;y\in{}^J\WW\},\{T_{y\i}\i;y\in{}^J\WW\}$ are two $H_{J,n}$-bases of $H_n$,
we have
$$T_y=\sum_{y'}c_{y,y'}T_{y'{}\i}\i,\qua T_{y\i}\i=\sum_{y'}d_{y,y'}T_{y'}$$
where $y,y'$ run over ${}^J\WW$ and $c_{y,y'},d_{y,y'}\in H_{J,n}$. For any $y$ we have
$$\align&\sum_{y'}\ov{a_{y,y'}}T_{y'{}\i}\i=v^{-2r}T_{y\i}\i C^\ss_{\uD\l}
=v^{-2r}\sum_{y''}d_{y,y''}T_{y''}C^\ss_{\uD\l}\\&
=v^{-2r}\sum_{y''}d_{y,y''}\sum_{y_1}a_{y'',y_1}T_{y_1}=v^{-2r}
\sum_{y''}d_{y,y''}\sum_{y_1}a_{y'',y_1}\sum_{y'}c_{y_1,y'}T_{y'{}\i}\i\endalign$$
hence 
$$\ov{a_{y,y'}}=v^{-2r}\sum_{y'',y_1}d_{y,y''}a_{y'',y_1}c_{y_1,y'}$$
for any $y,y'$. Hence it suffices to show that
$$\sum_y\z^{A'}(a_{y,\e(y)}[D])=
\sum_y\z^{A'}(\sum_{y'',y_1}d_{y,y''}a_{y'',y_1}c_{y_1,\e(y)}[D]).$$
By Lemma 31.8 (for $G',A'$ instead of $G,A$) we have
$$\z^{A'}(d_{y,y''}a_{y'',y_1}c_{y_1,\e(y)}[D])
=\z^{A'}(a_{y'',y_1}c_{y_1,\e(y)}[D]d_{y,y''})$$
Hence it suffices to show that
$$\sum_y\z^{A'}(a_{y,\e(y)}[D])=
\z^{A'}(\sum_{y,y'',y_1}a_{y'',y_1}c_{y_1,\e(y)}[D]d_{y,y''}).$$
We have 
$$\align&\sum_{y'}d_{\e(y),\e(y')}T_{\e(y')}=T_{\e(y)\i}\i=[D]T_{y\i}\i[D]\i\\&
=\sum_{y'}[D]d_{y,y'}[D]\i[D]T_{y'}[D]\i=\sum_{y'}[D]d_{y,y'}[D]\i T_{\e(y')}\endalign
$$
hence $d_{\e(y),\e(y')}=[D]d_{y,y'}[D]\i$. Hence it suffices to show that
$$\sum_y\z^{A'}(a_{y,\e(y)}[D])=\z^{A'}
(\sum_{y,y'',y_1}a_{y'',y_1}c_{y_1,\e(y)}d_{\e(y),\e(y'')}[D]).\tag b$$
From the definitions, $\sum_yc_{y_1,\e(y)}d_{\e(y),\e(y'')}$ is $1$ if $y_1=\e(y'')$ 
and is $0$, otherwise. Hence (b) holds. The lemma is proved.

\proclaim{Theorem 31.14}Let $D,P,L,G',D'$ be as in 29.1. Let $A$ be a character sheaf 
on $D$. Then $\res_D^{D'}A$ is a direct sum of character sheaves on $D'$.
\endproclaim
We can find $\cl\in\fs(\TT)$ and a sequence $\ss=(s_1,s_2,\do,s_r)$ in $\II$ such that
$s_1s_2\do s_r\uD\in\WW^\bul_\cl$ and such that $A$ is a direct summand of ${}^pH^i(K)$
for some $i\in\ZZ$, where $K=\bK^{\ss,\cl}_D[m]$, $m=r+\dim G$. Let 
$K'=\res_D^{D'}(K)$. For any $i$, let 
$$K_i={}^pH^i(K),K'_i=\res_D^{D'}(K_i).$$
For any character sheaf $A'$ on $D'$, let $b_{i,j}=(A':{}^pH^j(K'_i))$,
$b_j=(A':{}^pH^j(K'))$. From 28.12(b) we have
$${}^pH^j(K')={}^pH^j(\op_i\res_D^{D'}(K_i)[-i])=\op_i{}^pH^{j-i}(K'_i)$$
hence $b_j=\sum_ib_{i,j-i}$. Using 31.13, which is applicable since $\cl\in\fs_n$ for 
some $n\in\NN_{\kk}^*$, we get $b_j=b_{-j}$ for all $j$ hence
$$0=\sum_jjb_j=\sum_{i,j}jb_{i,j-i}=\sum_{i,j}(i+j)b_{i,j}.\tag a$$
From 28.17(a) we have $K_i=K_{-i}$. It follows that $b_{i,j}=b_{-i,j}$ so that
$\sum_{i,j}ib_{i,j}=0$. Introducing this into (a) we find $\sum_{i,j}jb_{i,j}=0$. From
28.12(b) and 30.6(b) we see that $b_{i,j}=0$ for all $j>0$. Therefore we have
$\sum_{i,j;j\le 0}jb_{i,j}=0$. Since $jb_{i,j}\le 0$ for all terms of the previous sum,
we must have $jb_{i,j}=0$ for all $i,j$. It follows that $b_{i,j}=0$ for $j\ne 0$. 
Since ${}^pH^j(K'_i)$ is a direct sum of character sheaves, by 29.15, it follows that
${}^pH^j(K'_i)=0$ for $j\ne 0$. In other words, for any $i$, $K'_i$ is a perverse sheaf
on $D'$ which is a direct sum of character sheaves. Since $A$ is a direct summand of 
$K_i$ for some $i$, we see that $\res_D^{D'}A$ is a direct summand of 
$\res_D^{D'}(K_i)=K'_i$ hence $\res_D^{D'}A$ is a perverse sheaf on $D'$ which is a
direct sum of character sheaves. The theorem is proved.

\proclaim{Corollary 31.15} Let $A$ be a character sheaf on $D$. Then $A$ is cuspidal 
(see 23.3) if and only if it is strongly cuspidal (see 23.3).
\endproclaim
We may assume that $A$ is cuspidal. Let $P,L,D'$ be as in 31.14 such that $P\ne G^0$. 
Since $A$ is cuspidal we have ${}^pH^i(\res_D^{D'}A)=0$ for all $i\ge 0$. By 31.14 we 
have ${}^pH^i(\res_D^{D'}A)=0$ for all $i\ne 0$. Hence ${}^pH^i(\res_D^{D'}A)=0$ for 
all $i$. It follows that $\res_D^{D'}A=0$. Thus, $A$ is strongly cuspidal. The 
corollary is proved.

\enddocument